\documentclass[11pt,leqno,twoside]{article}

 \usepackage{amsmath,amssymb,fancyhdr,makecmds,parskip,bm}%
\usepackage{fancyhdr}

\usepackage{fouriernc}

\def\ha{\frac{1}{2}} 
\def\thn{\hskip1.2pt}

\newcounter{dpage}

\def\twoupa4{%
\usepackage[landscape]{geometry}
\usepackage{2up}
\setlength{\paperwidth}{11.692in}
\setlength{\paperheight}{8.267in}
\setlength{\topmargin}{-0.75in}
\setlength{\headsep}{0.15in}
\setlength{\textheight}{7.3in}
\setlength{\oddsidemargin}{-0.75in}
\setlength{\evensidemargin}{-0.75in}
\setlength{\textwidth}{5.0in}
\source{\magstep0}{5.5in}{8.267in}
\target{\magstep0}{11.692in}{8.267in}
\def\oddno##1{\setcounter{dpage}{##1}\multiply \value{dpage}  by 2 \addtocounter{dpage}{-1}\thedpage}
\def\evenno##1{\setcounter{dpage}{##1}\multiply \value{dpage}  by 2 \addtocounter{dpage}{0}\thedpage}}

\newcounter{sequation}[section] 
\renewcommand{\thesequation}{\arabic{section}.\arabic{sequation}} 
 \newcounter{pequation} 
\renewcommand{\thepequation}{(\arabic{pequation})}%
\newcounter{problems} 
 \def\tl#1{\refstepcounter{sequation}\label{#1}{\text{\bm{$\thesequation$}}\thn}} 
  \def\dl#1{\refstepcounter{sequation}\label{#1}\leqno{\thesequation}} 
 \def\pdl#1{\refstepcounter{pequation}\label{#1}\leqno{\thepequation}} 
\def\tg#1{\refstepcounter{sequation}\label{#1}\tag*{$\thesequation$}} 
 \def\ptg#1{\refstepcounter{pequation}\label{#1}\tag*{$\thepequation$}}     

\newcommand{\prelemskip}{\vskip 5pt plus 1.5pt minus 1.5pt}
\newcommand{\postlemskip}{\vskip 0pt plus 1pt minus 0.5pt} 
\newenvironment{state}[1]{\par\prelemskip\noindent {\bf#1} %
\textit\bgroup\abovedisplayskip=5pt plus 3pt minus 1pt \belowdisplayskip=5pt plus 3pt minus 1pt}%
{\egroup \par\ifdim\lastskip<\smallskipamount \removelastskip\penalty55\postlemskip\fi}  
\newenvironment{proof}{\setcounter{pequation}{0}}{\,\,\nolinebreak$\Box$} 

 \renewcommand{\S}{\text{$\mathsection$}\thn}
\renewcommand{\epsilon}{\varepsilon}
\newcommand{\op}[1]{\operatorname{\text{\rm #1}}} 
\def\R{\mathbb{R}}      
\def\C{\mathbb{C}}   
\def\chir{\raise2pt\hbox{$\chi$}}
\def\ab#1#2#3{{\abovedisplayskip#1\belowdisplayskip#2
#3}} \DeclareMathSymbol{\pri}{\mathcal}{symbols}{"30}
\renewcommand{\prime}{\hskip0.8pt\pri\hskip-0.4pt}
\def\fr #1/#2{{\textstyle\frac{#1}{#2}}}

\def\thn{\hskip1.2pt}
\def\sing{\op{sing}}

\def\graph{\op{graph}}
\def\dist{\op{dist}}
\def\tint{\textstyle\int}
\def\res{\hbox{ {\vrule height .22cm}{\leaders\hrule\hskip.2cm} } }        
\def\sres{\hskip-0.7pt\raise-0.5pt\hbox{ {\vrule height .18cm}{\leaders\hrule\hskip.17cm} }\hskip-0.7pt }        
\def\reg{\op{reg}}
\def\spt{\op{spt}}
 \def\tsum #1{{\textstyle\sum} #1}

\fancypagestyle{plain}{%

}

\title {\vskip-.3in          The Symmetric Minimal Surface Equation
\vskip-.1in
\author{\scshape \normalsize     K.\ Fouladgar    \&    L.\ Simon\,\thanks{The research of the second author was partly
supported  by NSF DMS--0104049 \& DMS-0406209 at Stanford University; the work described here also includes part
of the Stanford doctoral thesis of the first author}}                                                  \date{\vspace{-.55in}}}%

\begin{document}

\maketitle  

\bigskip\bigskip

\thispagestyle{empty}

\bigskip\bigskip

\section{Introduction} \label{intro}

The Symmetric Minimal Surface Equation (SME) is the equation  
$$%
\mathcal{M}(u) = \frac{m-1}{u\sqrt{1+|Du|^{2}}} %
\dl{sme}  
$$%
on an open set $\Omega\subset\R^{n}$, where 
$$%
\mathcal{M}(u)=\tsum_{i=1}^{n}D_{i}\Bigl(\frac{D_{i}u}{\sqrt{1+|Du|^{2}}}\Bigr)  %
$$%
is the mean curvature operator on $\R^{n}$, $m\ge 2$ is an integer, and $u >0$ is of class $C^2$. Notice that this
equation is geometrically scale invariant: that is if $G(u)\subset\Omega\times(0,\infty)$ is the graph
$$%
G(u)=\{(x,u(x)):x\in\Omega\} %
\dl{graph} 
$$%
of a solution $u$ on $\Omega$ then, for each $\lambda>0$, $\lambda G(u)(=\{\lambda(x,u(x)):x\in\Omega\})$ is also the
graph of a solution; indeed it is the graph of the solution $u_{\lambda}(x)=\lambda u(x/\lambda)$ on the domain
$\lambda\Omega$. Analytically this is just the statement
$$%
u\text{ satisfies \ref{sme} on $\Omega\subset \R^{n}$}  \iff \lambda u(x/\lambda) \text{ satisfies \ref{sme} on  %
  $\lambda\Omega$ for each $\lambda>0$}.  %
\dl{gsi}
$$%
In fact geometrically the equation~\ref{sme} expresses the fact that the graph $G(u)$ of $u$ is a hypersurface with mean
curvature $H=\frac{m-1}{u\sqrt{1+|Du|^{2}}}$ at each of its points $(x, u(x))$.

The chief motivation here for the study of~\ref{sme} is the fact that if $u$ is a positive $C^{2}$ solution
of~\ref{sme} and if $S(u)$ is the ``symmetric graph'' of $u$, defined by
$$%
S(u) = \{(x,\xi)\in \Omega\times\R^{m}: |\xi|=u(x)\}, %
\dl{sym-graph}
$$%
then $S(u)$ is a minimal hypersurface (i.e.\ has zero mean curvature) in $\R^{n+m}$.

This is easily checked by observing that~\ref{sme} is the Euler-Lagrange equation for the functional
$$%
\mathcal{A}(u)= %
\sigma_{m-1}\tint_{\Omega}\sqrt{1+|Du|^{2}}\, u^{m-1}\,dx\,\,\, (\sigma_{m-1}={\cal{}H}^{m-1}(S^{m-1}))%
$$%
and, geometrically, $\mathcal{A}(u)$ represents the area functional for $S(u)$; that is, ${\cal{}A}(u)$ is the
$(n+m-1)$-dimensional Hausdorff measure ${\cal{}H}^{n+m-1}(S(u))$. This is clear because the integrand
$\sqrt{1+|Du|^{2}}\,u^{m-1}$ for ${\cal{}A}(u)$ is the Jacobian of the map $(x,\omega)\in \Omega\times S^{m-1}\mapsto
(x,u(x)\omega)\in \Omega\times \R^{m}$, and this map is a local coordinate representation for the symmetric graph
$S(u)$.  Since~\ref{sme} \ref{sme} expresses the fact that $u$ is stationary with respect to $\mathcal{A}$, we see that
$S(u)$ is stationary with respect to smooth symmetric deformations, and hence stationary with respect to all deformations
by a well-known principle (see e.g.\ \cite{Law72}).  (The latter principle here is just the natural generalization of the fact
that if a smooth hypersurface $\Sigma$ is rotationally symmetric about an axis and if $\Sigma$ is stationary with respect to
smooth rotationally symmetric compactly supported perturbations, then $\Sigma$ is minimal---i.e.\ stationary with
respect to all smooth compactly supported perturbations whether symmetric or not.)  Thus the smooth submanifold
$S(u)$ is stationary as a multiplicity~1 varifold in $\Omega\times(\R^{m}\setminus\{0\})$ and hence is a smooth
minimal submanifold of $\Omega\times(\R^{m}\setminus\{0\})$ as claimed.

\medskip

We say that a non-negative $C^{0}(\Omega)$ function is a \emph{singular solution} of the SME if it is locally the uniform
limit of $C^{2}$ positive solutions of the SME. More precisely:  

\smallskip
 
\tl{SMEsing} {\bf{}Definition.} $u:\Omega\to [0,\infty)$ is a singular solution of the SME in $\Omega$ if $u$ is
continuous on $\Omega$, if $\{x\in \Omega:u(x)=0\}\neq \emptyset$, and if 
$u=\lim_{j\to\infty}u_{j}$, uniformly on each compact subset of $\Omega$, where each $u_{j}$ is a positive
$C^{2}(\Omega)$ solution of~\ref{sme}.

\smallskip

If $u$ is  a regular or singular solution  of the SME in $\Omega$ then $\sing u$, \emph{the singular set} of $u$, is
defined to be $u^{-1}\{0\}$; that is,
$$%
\sing u= \{x\in\Omega: u(x)=0\},  %
\dl{sing-u}
$$%
Since by definition singular solutions are continuous at all points of $\Omega$, we see that $\sing u$ is closed as a subset
of $\Omega$, and of course $\sing u=\emptyset$ in case $u>0$ everywhere on $\Omega$.

\medskip 

\tl{SMEsing-rem} {\bf{}Remark:} For any regular or singular solution $u$ of \ref{sme} the symmetric graph $S(u)$ (as
in~\ref{sym-graph}) is a stationary multiplicity~1 varifold in the cylinder $\Omega\times\R^{m}$ and the singular set of
$S(u)$ indeed coincides with $\sing u\times \{0\}$; that is if $G$ is the graph of a singular solution $u$ of~\ref{sme} on
$\Omega$, and if $\sing S(u)$ is (as usual for stationary varifolds) defined to be the set of points $z\in\overline{S(u)}$ such
that there is no $\sigma>0$ such that $\overline{S(u)}\cap B_{\sigma}(z)$ is an $(n+m-1)$-dimensional embedded $C^{1}$
submanifold of $\Omega\times\R^{m}$, then
$$%
\sing S(u) =\sing u\times\{0\}\,(= \{(x,\xi)\in \Omega\times \R^{m}:\xi=0\text{ and }u(x)=0\}).
$$%
We will check this in Corollary~\ref{nulsing} after we have established the necessary preliminary area bounds.

\bigskip

It is not quite clear a-priori (but nevertheless true) that singular solutions $u$ of~\ref{sme} are in fact weak solutions
of~\ref{sme}, i.e.\ if $u$ is a singular solution of \ref{sme} then $u\in W^{1,1}_{\text{loc}}(\Omega)$, $1/u\in
L^{1}_{\text{loc}}(\Omega)$, and
$$%
\int_{\Omega}\Bigl(\sum_{i=1}^{n}\frac{D_{i}uD_{i}\zeta}{\sqrt{1+|Du|^{2}}}+ %
\frac{(m-1)\zeta}{u\sqrt{1+|Du|^{2}}}\Bigr)=0 \mbox{ for each } \zeta\in C^{1}_{c}(\Omega). %
\leqno{\ref{sme}^{\prime}}
$$%
We shall not explicitly use this fact here (although we will use~\ref{sme}$^{\prime}$ for regular solutions $u$), but the
interested reader can check that singular solutions $u$  also satisfy~\ref{sme}$^{\prime}$ by using the main regularity
theorem of \S\ref{dimsing}. (See Remark~\ref{weak-check} following Theorem~\ref{sing0}.)

\smallskip

As far as the \textit{existence} of singular solutions is concerned, we first note that there are no singular solutions in case
$n=1$: Indeed when $n=1$ there is a single variable $x$ and any positive solution $u$ on an interval $(a,b)$ satisfies the
ODE $\frac{u''}{1+(u')^{2}}=\frac{m-1}{u}$, so $u$ is strictly convex on $(a,b)$, and (after multiplying by $u'$ and
integrating) we see that $u^{1-m}(1+(u')^{2})^{1/2}=C$ (a positive constant). Thus $u$ is bounded below (indeed
$u^{m-1}\ge C^{-1}$). By the uniqueness and extension theorems for ODE's and the convexity of $u$ we then easily check
that, modulo a translation of the independent variable $x$, any positive solution extends to a maximal interval $(-d,d)$
with $0<d\le \infty$, where $u(x)=u(-x)$, $u$ takes its unique minimum at $x=0$ and $u(x)\to \infty$ as $|x|\to d$. 
Finally since the ODE is geometrically scale invariant it follows (again using the uniqueness theorem) that, after a translation
of the independent variable $x$, all solutions $u$ are just geometric rescalings $u(x)=\lambda\varphi(x/\lambda)$ for some
$\lambda>0$, where $\varphi$ is the unique (maximally extended) solution of the ODE with $\varphi(0)=1$,
$\varphi'(0)=0$ with $\varphi$ defined over some interval $(-d_{0},d_{0})$ where $0<d_{0}\le \infty$.  Now it is evident
that there can be no singular solution $u$, since otherwise there would be a sequence
$\lambda_{j}\varphi((x-x_{j})/\lambda_{j})$ with $\lambda_{j}\downarrow 0$ and $x_{j}\to 0$ which converges uniformly
on some open interval of $\R$.  This is impossible because if $d_{0}<\infty$ then $\varphi(x/\lambda_{j})$ is defined only
over the interval $(-d_{0}\lambda_{j},d_{0}\lambda_{j})\to \{0\}$, while if $d_{0}=\infty$ then
$\lambda_{j}\varphi(x/\lambda_{j})$ has derivative $\varphi'(x/\lambda_{j})$ which tends to $+\infty$ for $x>0$ and
$-\infty$ for $x<0$.

\smallskip

On the other hand in case $n\ge 2$ it is easy to give examples of singular solutions.  For instance one sees by direct
computation that ${u(x)\equiv(\frac{m-1}{n-1})^{\mspace{1mu}1/\mspace{1mu}2}} |x|$ is a solution of of the SME on
$\R^{n}\setminus\{0\}$, and, with a little more effort using ODE theory, it is straightforward to show that this $u$ is
locally the uniform limit of positive solutions in a neighborhood of $0$, and so is a singular solution in the sense introduced
above.  Notice that for this example the symmetric graph $S(u)$ is just the minimal cone
${|\xi|=(\frac{m-1}{n-1})^{1/\mspace{1mu}2}}|x|$ or in other words the cone $(n-1)|\xi|^{2}=(m-1)|x|^{2}$.  If $m=n=4$
this is the $7$-dimensional ``Simons cone'' over $S^{3}\times S^{3}$ which was the first example known of a singular area
minimizing hypersurface (\cite{BomDG69}, \cite{Sis68}).

\smallskip

We show in \S\ref{existence}, via a Leray-Schauder argument, that in fact there is a very rich class of singular solutions for
$n,m\ge 2$. In addition, one can quite easily modify the method of \cite{CafHS84} to show that the special example
$u_{0}(x)=(\frac{m-1}{n-1})^{1/\mspace{1mu}2}\,|x|$ generates a rich class of examples with isolated singular points
at~$0$, each of which is asymptotic to $u_{0}$ on approach to the singular point $0$.

Regular solutions of equation~\ref{sme} have been studied in \cite{DieH90}, \cite{DieH96}, with a different motivation
than the present one---the emphasis here is on the study of singular solutions of~\ref{sme} and the corresponding singular
minimal surfaces obtained from the symmetric graphs of solutions. Our aim is to develop the basic theory of such solutions,
showing on the one hand that there is a convenient general theory governing their qualitative behavior, and that on the
other hand, as mentioned above, there is a rich class of singular solutions, suggesting that this setting could be valuable in
improving our understanding of singular behavior of minimal submanifolds.

\smallskip

The main regularity results here (see \S\ref{grad-ests} and \S\ref{dimsing} below) are that singular solutions $u$
of~\ref{sme} are always locally Lipschitz in $\Omega$ with $\sing u$ having Hausdorff dimension less than or equal
$n-2$.

\section{Volume Bounds}\label{vol-bounds}

Here and subsequently, for $R>0$ and $y\in\R^{N}$,
$$%
B^{N}_{R}(y) = \{x\in\R^{N}:|x-y|<R\} \text{ (possibly abbreviated  $B_{R}(y)$ if $N$ is evident).}
$$%
We first want to present a lemma describing the basic volume bounds available for solutions of the SME (i.e.\ solutions
of~\ref{sme}). 

As a preliminary to this, recall (see e.g.\ \cite{Sim83}) that if $U\subset \R^{n+1}$ is open and $\Sigma$ is a smooth
embedded hypersurface (for the moment with no singular set in $U$, i.e.\ $(\overline \Sigma\setminus \Sigma)\cap
U=\emptyset$) then
$$%
e^{\Lambda \rho} \rho^{-n}{\cal{}H}^{n}(\Sigma\cap B_{\rho}(y))\text{ is increasing in $\rho$, $\rho\in (0,R]$}, %
\dl{monotonecor}
$$%
provided $\overline{\!B}_{R}(y)\subset U$ and $\sup_{B_{R}(y)}|H|\le \Lambda$, where $H$ is the mean curvature of
$\Sigma$. In particular
$$%
\mu(\Sigma\cap B_{\rho}(y))\ge e^{-\Lambda \rho}\omega_{n}\rho^{n},\quad \rho\in (0,R], \, y\in\spt \mu. 
\dl{monotonecor1}
$$%

The main volume bounds for solutions of~\ref{sme} are as follows:
   
\begin{state}{\tl{vol-bds}\,Lemma.}%
  There is a constant $C=C(m,n)$ such that if $u$ is a regular or singular solution of \emph{\ref{sme}} on a domain
$\Omega\subset \R^{n}$ and if $y=(x_{0},u(x_{0}))\in G= \graph u$ and $\overline{\!B}_{\rho}(y)\cap G$ is compact (so
the boundary of $G$ does not intersect $B_{\rho}(y)$), then
$$%
C^{-1}\rho^{n}\le \tint_{S_{\!\rho/2}(y)}\sqrt{1+|Du|^{2}} \,dx\le C \rho^{n} %
\leqno{\emph{(i)}}
$$%
and
$$%
\tint_{S_{\!\rho/2}(y)}\sqrt{1+|Du|^{2}}\,u^{m-1}\,dx %
\ge C^{-1}(\rho+u(x_{0}))^{m-1}\rho^{n}, %
\leqno{\emph{(ii)}}
$$%
where we use the notation $S_{\sigma}(y)=\{x:(x,u(x))\in B_{\sigma}^{n+1}(y)\cap G\}$; that is,
$$%
S_{\sigma}(y)=\{x\in \Omega:\sqrt{|x-x_{0}|^{2}+ (u(x)-u(x_{0}))^{2}} < \sigma\}.
$$%
\end{state}
 
{\bf \tl{volbds-rem} Remark:} Observe that the above lemma implies that there is
$\varepsilon_{0}=\varepsilon_{0}(m,n)>0$ such that if $u$ is a regular or singular solution of~\ref{sme} on a ball
$B_{\rho}(x_{0})\subset\R^{n}$ then there is at least one point $x$ in $B_{\rho/2}(x_{0})$ with
$u(x)>\varepsilon_{0}\rho$, because otherwise by the right inequality of~(i) we would have
\ab{3pt}{3pt}{%
$$%
\tint_{S_{\rho}(x_{0},u(x_{0}))}\sqrt{1+|Du|^{2}}\,u^{m-1}\,dx\le C\varepsilon_{0}^{m-1} \rho^{n+m-1},
$$}%
contradicting the bound (i\hskip-.05pti) of the lemma.

\smallskip

We can now check the claim made in Remark~\ref{SMEsing-rem}:
\begin{state}{\tl{nulsing}\,Corollary.}%
  If $u$ is a regular or singular solution of~\emph{\ref{sme}} in an open $\Omega\subset\R^{n}$, then the symmetric graph
$S(u)$ (as in~\emph{\ref{sym-graph}}) is a multiplicity~$1$ stationary $(n+m-1)$-dimensional varifold in
$\Omega\times\R^{m}$ with $\sing S(u)=\sing u\times\{0\}=\overline{S(u)}\setminus S(u)$, where $\overline{S(u)}$ is the
closure of $S(u)$ in $\Omega\times\R^{m}$.  Furthermore in both the regular and singular case $u$ is $C^{\infty}$ on
the (open) set of points where it is positive, and $G(u)\cap (\Omega\times(0,\infty)), S(u)\cap
(\Omega\times(\R^{m}\setminus \{0\}))$ are $C^{\infty}$ hypersurfaces in $\R^{n+1}$ and $\R^{n+m}$ respectively.
  \end{state}

\smallskip

  \begin{proof}{\bf{}Proof:} In case $u>0$ in $\Omega$ (i.e.\ the case when $\sing
u=\emptyset$),  $S(u)$ is a smooth $(n+m-1)$-dimensional minimal submanifold of
$\Omega\times(\R^{m}\setminus\{0\})$,  as discussed in \S\ref{intro}. Hence if $u>0$ in  $\Omega$ then $S(u)$ is
stationary, as a multiplicity~$1$ varifold, in $\Omega\times\R^{m}$ in accordance with the discussion of \S\ref{intro}.

If now $u$ is a singular solution of the SME on $\Omega$, then there is a sequence $u_{j}$ of positive solutions
of~\ref{sme}  with $u_{j}\to u$ uniformly on compact subsets of $\Omega$.  Using the gradient estimates of
\cite{Sim76} (applied in the case when the functions $A_{i},B$ of \cite{Sim76} satisfy
$A_{i}(x,u,Du)=D_{i}u/\sqrt{1+|Du|^{2}},\,|B(x,u,Du)|\le C/\sqrt{1+|Du|^{2}}$) together with quasilinear elliptic regularity
theory \cite[Chapter 10]{GT}, we deduce that, for each $k$, $u$ is the $C^{k}$ limit of $u_{j}$ in a neighborhood of each
point where $u$ is non-zero.  So indeed $G(u)\cap (\Omega\times(0,\infty)), S(u)\cap (\Omega\times(\R^{m}\setminus
\{0\}))$ are $C^{\infty}$ as claimed.

Since we have the local area bounds of Lemma~1 and since each $S(u_{j})$ is stationary as a multiplicity~1 varifold in
$\Omega\times\R^{m}$,  by the Allard compactness theorem for integer multiplicity varifolds there is a subsequence
and a limiting integer multiplicity varifold $V$ which is stationary in $\Omega\times\R^{m}$, and by the above
discussion $V$ is smooth in $\Omega\times(\R^{m}\setminus\{0\})$.  Furthermore if $B^{n}_{\rho}(x_{0})\subset\Omega$
then, for each $\sigma\in (0,\rho/2)$, $B^{n+m}_{\rho/2}(x_{0},0)\cap \{(x,y)\in\R^{n}\times\R^{m}:|y|\le \sigma/2\}$ can
be covered by balls $B_{\sigma}(y_{1},0),\ldots, B_{\sigma}(y_{N},0)$ with $N\le C (\rho/\sigma)^{n}$, and the upper
bound of Lemma~\ref{vol-bds}(i) gives ${\cal{}H}^{n+m-1}(S(u_{j})\cap B_{\sigma}(y_{k}))\le C\sigma^{n+m-1}$ for each
$k$, so ${\cal{}H}^{n+m-1}(S(u_{j})\cap \{(x,y)\in\R^{n}\times\R^{m}:|y|\le \sigma/2\}\le C\rho^{n}\sigma^{m-1}\to 0$ as
$\sigma\downarrow0$.  Since by the discussion above we also know the $S(u_{j})$ converges to $S(u)$ in the $C^{k}$
sense on $\{(x,y)\in\R^{n}\times\R^{m}:|y|> \sigma/2\}$ for each $\sigma>0$, we thus have that $S(u)$ is the varifold
limit of $S(u_{j})$ (assuming $S(u)$ is viewed as a multiplicity~$1$ varifold in $\R^{n+m}$), so indeed $S(u)=V$ and
$S(u)$ is stationary as claimed.

To complete the proof we have to check $\sing S(u)=u^{-1}(0)\times\{0\}$. Certainly $\sing S(u)\subset
u^{-1}\{0\}\times\{0\}$ because, as already mentioned, $S(u)$ is smooth in $\Omega\times(\R^{m}\setminus\{0\})$. So we
have only to check that no point $(x,0)\in u^{-1}\{0\}\times\{0\}$ is a regular point of $S(u)$ and without loss of generality
we just check the case $x=0$. If $0$ is a regular point of $S(u)$ then there is $\sigma>0$ with $\Sigma=S(u)\cap
B^{n+m}_{\sigma}(0)$ a smooth $(n+m-1)$-dimensional embedded manifold of $\R^{n}\times\R^{m}$.  Evidently, since
$\Sigma$ is invariant under orthogonal transformations of the last $m$ variables, $T_{w}\Sigma +
(\R^{n+1}\times\{0\})=\R^{n+m}$ for each $w\in \Sigma\cap(\R^{n+1}\times\{0\})\setminus(\R^{n}\times\{0\})$, and by
\ref{volbds-rem} each $w\in \Sigma\cap (\R^{n}\times\{0\})$ is the limit of a sequence $w_{j}\in
\Sigma\cap(\R^{n+1}\times\{0\})$, so in fact $T_{w}\Sigma + (\R^{n+1}\times\{0\})=\R^{n+m}$ for each $w\in
\Sigma\cap(\R^{n+1}\times\{0\})$.  Therefore by transversality theory $\Sigma$ intersects $\R^{n+1}\times\{0\}$
transversely with intersection a smooth embedded $n$-dimensional submanifold, but clearly
$\Sigma\cap(\R^{n+1}\times\{0\})=(G(u)\cup G(-u))\cap B^{n+1}_{\sigma}(0)$, so $(G(u)\cup G(-u))\cap
B^{n+1}_{\sigma}(0)$ is a smooth embedded $n$-dimensional submanifold of $\R^{n+1}$.  By
Remark\thn\ref{volbds-rem} the tangent space of this submanifold at $0$ is distinct from $\R^{n}\times\{0\}$ and hence,
again by transversality theory and taking a smaller $\sigma$ if necessary, we see that $B^{n+1}_{\sigma}(0)\cap (G(u)\cup
G(-u))\cap(\R^{n}\times\{0\})$ is a smooth connected embedded $(n-1)$-dimensional manifold $\Gamma$ and
$B_{\sigma}^{n+1}(0)\cap G(u)\setminus \Gamma$ is connected.  But $B^{n+1}_{\sigma}(0)\cap G(u)\setminus \Gamma$
can be written as the disjoint union $(B^{n+1}_{\sigma}(0)\cap\graph(u|U_{+}))\cup(B^{n+1}_{\sigma}(0)\cap\graph
(u|U_{-}))$, where $U_{\pm}$ are the two connected components of $B^{n}_{\sigma}(0)\setminus\Gamma$, and, by
continuity of $u$ at $0$, both sets in this union are non-empty, contradicting the connectedness of
$B^{n+1}_{\sigma}(0)\cap G(u)\setminus \Gamma$.~\end{proof}

\medskip

\begin{proof}{\bf{}Proof of Lemma \ref{vol-bds}:} First we prove the upper bound in~(i).  We can assume that $u$ is a
regular solution, otherwise apply the argument to an approximating sequence of regular solutions $u_{j}$.  Replace
$\zeta$ in the weak form \ref{sme}$^{\prime}$ of \ref{sme} by $\gamma(u)\zeta$, where $\zeta\in C^{1}_{c}(B_{\rho}(y))$
is non-negative and $\gamma(t)$ is the piecewise linear function $\R\to \R$ which is zero for $t\le (u(x_{0})-\rho/2)_{+}$,
slope $1$ for $t\in (u(x_{0})-\rho/2,u(x_{0})+\rho/2)$ and $\gamma(t)=\rho$ for $t\ge u(x_{0})+\rho/2$.  Then since
\ab{3pt}{3pt}{%
$$%
\tsum_{i=1}^{n}\tfrac{D_{i}u}{\sqrt{1+|Du|^{2}}}D_{i}u\ge \sqrt{1+|Du|^{2}}-1,
$$}%
we deduce that
\begin{align*}%
  &\tint_{\{x\in B_{\rho}(y):-\rho/2<u(x)-u(x_{0})<\rho/2\}}\zeta\sqrt{1+|Du|^{2}}\,dx \\ %
  &\hskip0.6in \le \tint_{B_{\rho}(y)}\bigl(m\zeta+\gamma(u) \tsum_{i=1}^{n}  %
\frac{D_{i}u}{\sqrt{1+|Du|^{2}}}D_{i}\zeta\bigr)\,dx \le C\rho^{n}+\rho\tint_{B_{\rho}(y)} |D\zeta|\,dx, %
\end{align*}%
and choosing $\zeta$ to be a standard cut-off function in $B_{\rho}(y)$ with $\zeta\equiv 1$ in $B_{\rho/2}(y)$ we then
obtain the desired bound.

Next we observe that since the symmetric graph $S(u)$ of $u$ is a minimal hypersurface, then the lower
bound \ref{monotonecor1} applies to give
$$%
{\cal{}H}^{n+m-1}(S(u)\cap Q(B^{n+m}_{\sigma}(y,u(x_{0}),0))\ge \omega_{m+n-1}\sigma^{m+n-1} %
\pdl{lower-bd}
$$%
for any orthogonal transformation $Q$ of $\R^{n+m}$ which acts as the identity on the first $n$-coordinates (i.e.\ the
$x$ coordinates) and any $\sigma\in (0,\rho)$, where $\omega_{n+m-1}$ is the measure of the unit ball in
$\R^{n+m-1}$. Also, using the symmetry of $S(u)$, we see that
$$%
\{(x,\xi):|x-y|<\rho/4,u(x_{0})-\rho/4<|\xi|<u(x_{0})+\rho/4\}(\subset S_{\rho/2}(x_{0}))
\pdl{inclusion}
$$%
contains pairwise disjoint balls $Q_{j}B^{n+m}_{\rho/4}(y,u(x_{0}),0)$, where $Q_{j}$, $j=1,\ldots,N$, are orthogonal
transformations of $\R^{n+m}$ which act as the identity in the first $n$ coordinates and
$$%
N\ge\max\{1, C(u(x_{0})+\rho)^{m-1}/\rho^{m-1}\},\quad C=C(n,m)>0.
$$%
We then have  by \ref{lower-bd} and \ref{inclusion}
\ab{3pt}{3pt}{%
$$%
  {\cal{}H}^{n+m-1}(S_{\rho/2}(x_{0}))  \ge C (u(x_{0})+\rho)^{m-1}\rho^{n},
$$}%
which is the required inequality~(ii).

Finally to prove the lower bound in~(i), we consider cases $u(x_{0})\ge \rho/4$, $u(x_{0})<\rho/4$.  If $u(x_{0})\ge
\rho/4$, the mean curvature $H$ of graph $G$ in $S_{\rho/8}(y)$ satisfies $\rho |H|\le (m-1)\rho/u\le 8(m-1)$, and a
standard consequence of the monotonicity inequality for surfaces of such bounded mean curvature is exactly that
\ab{4pt}{4pt}{%
$$%
e^{8(m-1)\sigma/\rho}\sigma^{-n}{\cal{}H}^{n}(G\cap S_{\sigma}(y))\text{ is increasing for }\sigma\in(0,\rho/8],
$$}%
and so the inequality $\tint_{S_{\rho/2}(y)}\sqrt{1+|Du|^{2}}\,dx\ge \tint_{S_{\rho/8}(y)}\sqrt{1+|Du|^{2}}\,dx\ge
C^{-1}\rho^{n}$ follows.  On the other hand if $u(x_{0})<\rho/4$ then $u<\rho/4+\rho/2=3\rho/4<\rho$ in $S_{\rho/2}$
and we can use the bound~(ii) to give
\ab{3pt}{3pt}{%
$$%
C\rho^{m-1}\tint_{S_{\!\rho/2}(y)}\sqrt{1+|Du|^{2}}\,dx\ge \tint_{S_{\!\rho/2}(y)}\sqrt{1+|Du|^{2}}\,u^{m-1}\,dx\ge
C^{-1}\rho^{m+n-1}, 
$$}%
and so again  $\tint_{S_{\rho/2}(y)}\sqrt{1+|Du|^{2}}\,dx\ge C^{-1}\rho^{n}$. Thus the proof of Lemma~\ref{vol-bds}
is complete.~\end{proof}

\section{H\"older Continuity}

Here we establish H\"older estimates for regular and singular solutions of the SME (equation~\ref{sme}). These will be
important in the proof of both the existence result in \S\ref{existence} and also in the proof of the gradient estimate in
\S\ref{grad-ests}.
 
\begin{state}{\tl{holder1} Theorem {\bf (H\"older Continuity).}} %
  Let $u\in C^{2}(\Omega)\cap C^{0}(\overline\Omega)$ satisfy~\emph{\ref{sme}} in $\Omega$ and $(u-\varphi)|\partial
\Omega=0$, where $\varphi$ is $C^{1,1}(\R^{ n})$ and $0\le\varphi\le u$ in $\Omega$.  Then $u$ is H\"older
continuous with exponent $\frac{1}{2}$ on $\overline \Omega$, and in fact
$$%
|u(x)-u(y)| \le C|x-y|^{\frac{1}{2}},\,\quad x,\,y\in \overline\Omega,
$$%
where $C=C(M,n)$; $C$ does not depend on $\Omega$. Also
$$%
\sup_{\Omega}|D(u-\varphi)^{2}| \le C
$$%
where again $C$ depends only on $M,n$.
\end{state}

Before giving the proof, we observe that the above theorem directly implies a local interior H\"older estimate for
solutions of~\ref{sme} and a local bound on the gradient of $u^{2}$:

\begin{state}{\tl{holdercor} Corollary.} %
  If $u$ is a regular or singular solution of~\emph{\ref{sme}} in $B_{\rho}(z)$ then $u$ is locally H\"older continuous with
exponent $\frac{1}{2}$ in $B_{\rho}(z)$, and 
\begin{align*}   
  |u(x)-u(y)| &\le C|x-y|^{\frac{1}{2}},\,\quad x,y\in B_{\rho/2}(z)\\  %
\sup_{B_{\rho/2}(z)}|Du^{2}| &\le C,
\end{align*} 
where $C$ depends only on $n,m$ and $M/\rho$, where $M$ is any upper bound for $u$ on $B_{\rho}(z)$.
\end{state}

{\bf{}Remark:} We note also that the above corollary implies that the gradient of $u$ is bounded on the subset of
$B_{\rho/2}(z)$ where $u(x)\ge \varepsilon$ for each $\varepsilon>0$.

\medskip

\begin{proof}{\bf Proof of the Corollary~\ref{holdercor}:} We can suppose that $u$ is a regular solution on $B_{\rho}(z)$
(otherwise uniformly approximate $u$ be positive solutions $u_{j}$ and apply the estimates to $u_{j}$). 

Choose $\psi$ to be any non-negative smooth function with $\psi\equiv0$ on $B_{\rho/2}(z)$, $\psi|\partial
B_{\rho}(z)=2\sup_{B_{\rho}(z)}u$ and
$$%
|D^{j}\psi|\le C(n)\rho^{-j}\sup_{B_{\rho}(z)}u \text{ for }j=1,2,
$$%
and then apply \ref{holder1} with $\widetilde \Omega$ in place of $\Omega$, where $\widetilde\Omega= \{x\in
B_{\rho}(z) : u(x) > \psi(x)\}$. The lemma gives $\sup_{B_{\rho/2}(z)}|Du^{2}|\le C$ (where $C$ depends only on
$\rho^{-1}M$ and $n$), which implies the H\"older estimate $|u(x)-u(y)| \le C|x-y|^{\frac{1}{2}},\, x,y\in B_{\rho/2}(z)$. 
\end{proof}

\medskip

\begin{proof}{\bf Proof of the Theorem \ref{holder1}:} We use a modification of the method of~\cite{KorS89}, which was
used to establish such H\"older estimates for solutions of $\mathcal{M}(u)=H(x,u)$ in the case when $H(x,z)$ is
\textit{increasing} in $z$, which is the correct sign for application of the maximum principle. The special form of the right
side, including in particular the factor $1/\sqrt{1+|Du|^{2}}$, allows the method of \cite{KorS89} to be successfully modified
to the present setting, even though $1/u$ is a decreasing function of $u$, as we now show.

We let $\Delta_{G}$ denote the Laplace-Beltrami operator on the graph $G=\{(x,u(x)) : x\in \Omega\}$ expressed in the
natural coordinates $x\in \Omega$: thus using the notation
$$%
v = (1+|Du|^{2})^{1/2},
$$%
we let $\displaystyle\nu=(\nu_{1} ,\dotsc, \nu_{n+1})=v^{-1}(-Du,1)$ be the upward pointing unit normal of $G(u)$ and
$\Delta_{G}\psi= v^{-1}\sum_{i,j=1}^{n} D_{i}(v(\delta_{ij}-\nu_{i}\nu_{j})D_{j}\psi)$.  Keeping in mind the classical
identity
$$%
v^{-1}\tsum_{i=1}^{n}D_{i}(v(\delta_{ij}-\nu_{i}\nu_{j}))=H\nu_{j}
$$%
for the mean curvature $H=-\sum_{i=1}^{n}D_{i}\nu_{i}(=(m-1)\frac{\nu_{n+1}}{u}\mbox{ by }\ref{sme}$), we see
that this can in fact alternatively be written
$$%
\Delta_{G}\psi=(\delta_{ij}-\nu_{i}\nu_{j})D_{i}D_{j}\psi + (m-1)\frac{\nu_{n+1}}{u}\nu_{j}D_{j}\psi,
$$%
where, here and subsequently, repeated indices are summed from $1$ to $n$. Using the abbreviations
$g^{ij}=\delta_{ij}-\nu_{i}\nu_{j}$, $u_{k}=D_{k}u, u_{ik}=D_{i}D_{k}u$, we also directly compute
\begin{align*}%
  \Delta_{G} \nu_{n+1} %
 &= v^{-1}D_{j}(vg^{ij}D_{i}\nu_{n+1}) %
  =-v^{-1}D_{j}\bigl(v^{-1}g^{ij}{u_{k}\over{}v}u_{ki}\bigr)\\ %
 &=-{u_{k}\over{}v^{2}}D_{j}(v^{-1}g^{ij}u_{ki}) - v^{-3}g^{ij}g^{k\ell}u_{ki}u_{\ell j} %
  &\hskip-0.6in\bigl(\text{as }D_{j}\bigl({u_{k}\over{}v}\bigr)=v^{-1}g^{k\ell}u_{\ell j}\bigr)\\ %
  &= -{u_{k}\over{}v^{2}}D_{j}(v^{-1}g^{ij}u_{ki}) - |A|^{2}\nu_{n+1} %
  &\hskip-1.1in(\text{where }|A|^{2}=v^{-2}g^{ij}g^{k\ell}u_{ik}u_{j\ell})\\ %
  &= -{u_{k}\over{}v^{2}}D_{j} %
  \bigl(D_{k}\bigl({u_{j}\over{}v}\bigr)\bigr) - |A|^{2}\nu_{n+1} %
  &\hskip-1.1in\bigl(\text{since } %
  D_{k}\bigl({u_{j}\over{}v}\bigr)=v^{-1}g^{ij}u_{k i}\bigr)\\ %
  &= -{u_{k}\over{}v^{2}}D_{k} %
  \bigl(D_{j}\bigl({u_{j}\over{}v}\bigr)\bigr) - |A|^{2}\nu_{n+1} \\ %
  &=-{u_{k}\over{}v^{2}}D_{k}\bigl({(m-1)\nu_{n+1}\over{}u}\bigr) - |A|^{2}\nu_{n+1} %
  &\hskip-1.1in \text{(by \ref{sme})} \\%
 &=\nu_{n+1}\nu_{k}D_{k}\bigl({(m-1)\nu_{n+1}\over{}u}\bigr) - |A|^{2}\nu_{n+1}. %
\end{align*}%
($|A|^{2}=v^{-2}g^{ij}g^{k\ell}u_{ik}u_{j\ell}$ is geometrically the squared length of the second fundamental form of
$G$).  Thus, in summary,
$$%
\Delta_{G}\nu_{n+1} + |A|^{2}\nu_{n+1} = (m-1) \nu_{n+1} \nu_{j}D_{j}(\nu_{n+1}/ u). %
\pdl{laplace}
$$%
Now we define
$$%
\eta = e^{K(u-\varphi)}-1  \quad (K>0 \text{ to be chosen})
$$%
and let $M=\sup {\frac{\eta}{\varepsilon+\nu_{n+1}}}$ (where $\varepsilon>0$ will be allowed to approach $0$ shortly)
and we observe that then $\eta - M(\varepsilon+\nu_{n+1})$ has a maximum value of $0$ which is attained at some point
$x_{0}\in \Omega$. Thus
$$%
D(\eta-M(\varepsilon+\nu_{n+1}))(x_{0})=0 \mbox{ and } %
\Delta_{G}(\eta-M(\varepsilon+\nu_{n+1}))(x_{0})\le 0. %
$$%
On the other hand we can directly compute $\Delta_{G}(\eta-M(\varepsilon+\nu_{n+1}))(x_{0})$.  In this
computation we let $h(u)=(m-1)/u$ (so the mean curvature $H=-D_{i}\nu_{i}$ of $G$ is just $h(u)\nu_{n+1}$),
$g^{ij}=(\delta_{ij}-\nu_{i}\nu_{j})$, subscripts (like $i,j$ in $\varphi_{i},\varphi_{ij}$) denote partial derivatives,
and we use the summation convention that repeated indices are summed from $1$ to $n$:
\begin{align*}%
  \ptg{key}  &\hskip-0.5in\Delta_{G}(\eta-M(\varepsilon+\nu_{n+1}))\\ %
 &= g^{ij}\eta_{ij}+h\nu_{n+1}\nu_{j}\eta_{j}-M\Delta_{G}\nu_{n+1} \\
&\ge g^{ij}\eta_{ij}+h\nu_{n+1}\nu_{j}\eta_{j} - M\nu_{n+1}\nu_{j}D_{j}(h\nu_{n+1}) %
  \hskip.2in\mbox{ (by \ref{laplace} above)}\\%
&= %
  g^{ij}\eta_{ij}+\nu_{n+1}\nu_{j}D_{j}(h(\eta-M\nu_{n+1})) -\nu_{n+1}\eta\nu_{j}D_{j}h\\
  \hspace{0.7in}&= K^{2}e^{K(u-\varphi)}g^{ij}(u_{i}-\varphi_{i})(u_{j}-\varphi_{j}) %
  +Ke^{K(u-\varphi)}g^{ij}(u_{ij}-\varphi_{ij}) \\%
 &\hskip 0.8in%
  -\nu_{n+1}\eta\nu_{j}D_{j}h +\nu_{n+1}\nu_{j}D_{j}\bigl(h(\eta-M\nu_{n+1})\bigr)\\%
 &= K^{2}e^{K(u-\varphi)}g^{ij}(u_{i}-\varphi_{i})(u_{j}-\varphi_{j}) %
  -Ke^{K(u-\varphi)}g^{ij}\varphi_{ij} \\ %
 &\hskip 0.4in +Ke^{K(u-\varphi)}h- \nu_{n+1}\eta\nu_{j}D_{j}h %
 +\nu_{n+1}\nu_{j}D_{j}\bigl(h(\eta-M\nu_{n+1})\bigr) %
 \end{align*}%
where we used the fact that $g^{ij}u_{ij}=h$ (by~\ref{sme}). Now at the point $x_{0}$ where
$\eta-M(\varepsilon+\nu_{n+1})$ has its maximum value of zero, we have $D_{j}(\eta-M\nu_{n+1})=0$ and
$\eta-M\nu_{n+1}=M\varepsilon$, so
$$%
\nu_{n+1}\nu_{j}D_{j}(h(\eta-M\nu_{n+1}))=M\varepsilon\nu_{n+1}\nu_{j}D_{j}h = %
M(m-1)\varepsilon|Du|^{2}u^{-2}v^{-2}\ge 0, %
$$%
and thus the crucial remaining point is in the sign of the term $Ke^{K(u-\varphi)}h -\nu_{n+1}\eta\nu_{j}D_{j}h$; since
$e^{Kt}-1\le Kte^{Kt}$ for $t\ge 0$ (which in particular guarantees $\eta\le K(u-\varphi)e^{K(u-\varphi)}\le
Kue^{K(u-\varphi)}$), and since $\nu_{n+1}\nu_{j}D_{j}h=(m-1)u^{-2}v^{-2}|Du|^{2}$, we see that in fact
$$%
Ke^{K(u-\varphi)}h-\nu_{n+1}\eta\nu_{j}D_{j}h\ge Khe^{K(u-\varphi)} %
\Bigl(1-\frac{|Du|^{2}}{v^{2}}\Bigr)=Khe^{K(u-\varphi)} v^{-2}\ge 0.%
$$%
Thus, at the point $x_{0}$ where $\eta-M(\epsilon+\nu_{n+1})$ takes its zero maximum value, \ref{key} gives
\begin{align*}   
  0 &\ge \Delta_{G}(\eta-M(\varepsilon+\nu_{n+1})) %
  \ge K^{2}e^{K(u-\varphi)}g^{ij}(u_{i}-\varphi_{i})(u_{j}-\varphi_{j}) %
  -Ke^{K(u-\varphi)}g^{ij}\varphi_{ij} \\%
  & = K^{2}e^{K(u-\varphi)}\Bigl(\frac{|Du|^{2}}{1+|Du|^{2}}- %
  2\frac{u_{j}\varphi_{j}}{1+|Du|^{2}} +g^{ij}\varphi_{i}\varphi_{j}\Bigr) %
  -Ke^{K(u-\varphi)}g^{ij}\varphi_{ij}\\ %
  & \ge Ke^{K(u-\varphi)}\Bigl(K \Bigl(\frac{|Du|^{2}}{1+|Du|^{2}}- %
    \gamma\frac{2|Du|}{1+|Du|^{2}}\Bigr)-\gamma\Bigr), %
\end{align*} 
where $\gamma=\sup_{\Omega}(|D\varphi|+\sum_{i,j}|\varphi_{ij}|)$, so we conclude
$$%
\frac{|Du|^{2}}{1+|Du|^{2}}-\frac{2\gamma|Du|}{1+|Du|^{2}} \le \gamma/K
$$%
at the point $x_{0}$ where $\eta/(\varepsilon+\nu_{n+1})$ has its maximum. Using Cauchy's inequality on the left we see
that then
$$%
\frac{\frac{1}{2}|Du(x_{0})|^{2}-2\gamma^{2}}{1+|Du(x_{0})|^{2}} \le \gamma/K,
$$%
and selecting $K=4\gamma$ we see that then
$$%
|Du(x_{0})|^{2}\le 8\gamma^{2}+1, 
$$%
hence $|Du(x_{0})|\le 4(\gamma +1)$. Thus, with the above choice $K=4\gamma$, we get
$M=\sup_{\Omega}\eta/(\varepsilon+\nu_{n+1}) \le \sup_{\Omega}e^{4\gamma(u-\varphi)}/(\varepsilon +
(4\gamma+4)^{-1})$, so that letting $\varepsilon\downarrow 0$ we have
$$%
\sup_{\Omega}\bigl(|Du|(e^{4\gamma(u-\varphi)}-1)\bigr) \le %
\sup_{\Omega}\bigl(\nu_{n+1}^{-1}(e^{4\gamma(u-\varphi)}-1)\bigr) \le %
4(\gamma+1)\sup_{\Omega}e^{4\gamma(u-\varphi)}.
$$%
The required inequality $|D(u-\varphi)^{2}|\le C$ now follows because $e^{4\gamma(u-\varphi)}-1\ge
4\gamma(u-\varphi)$.
\end{proof}

\section{Existence Results} \label{existence}

In this section $n\ge 2$ and we show there exists quite a rich class of singular solutions of the SME (i.e.\ \ref{sme}) on
any uniformly convex $C^{2,\alpha}$ domain $\Omega\subset\R^{n}$.

We start with the following definitions:

{\bf{}\tl{st-pos-non-solv} Definitions:} Boundary data $\varphi\in C^{0}(\partial\Omega)$ is said to be \emph{strongly
positive} if there exists $\eta >0$ such that $\inf_\Omega u\geq \eta$, whenever $u>0$ is a $C^{0}(\overline\Omega)\cap
C^{2}(\Omega)$ solution of \ref{sme} with $\min_{\partial\Omega}u\ge \varphi$, and $\varphi$ is said to be
\emph{non-solvable} if there exists no $C^{2}(\Omega)\cap C^{0}(\overline\Omega)$ solution of \ref{sme} with
$u|\partial\Omega=\varphi$.

Of course there exists such non-solvable data; indeed if $\epsilon>0$ and
$\sup_{\partial\Omega}\varphi<\epsilon$ then by the maximum principle any solution $u$ of \ref{sme} with boundary data
$\varphi$ would also satisfy $\sup_{\Omega} u\le \epsilon$, and for small enough $\epsilon=\epsilon(m,\Omega)>0$ this is
impossible by Remark\thn\ref{volbds-rem}.

There also exists strongly positive data for any given $C^{2}$ domain $\Omega\subset\R^{n}$ for $n\ge 2$---indeed there
is $K>0$ such that any $\varphi$ with $\min_{\partial\Omega}\varphi\ge K$ is strongly positive, as one easily checks as
follows:

First consider solutions of the SME \ref{sme} on $\{x\in \R^{n}:|x|>1\}$ and which are functions of $r=|x|$, so
$u_{0}(x)=\psi(r)$ with $\psi\in C^{\infty}(0,\infty)$.  One can check using ODE theory that such solutions exists
and indeed there is a unique such solution which satisfies
\ab{3pt}{3pt}{%
$$%
\psi(1+)\,\bigl(=\lim_{r\downarrow 1}\psi(r)\bigr)=0,\,  \psi'(1+)=+\infty, \text{ and } %
\lim_{r\to\infty}|\psi(r) - \sqrt{\tfrac{m-1}{n-1}}\,r|=0,  %
\dl{ext-ode-soln}
$$}%
With this solution $\psi$ and $\lambda>0$, let
$$%
\psi_{\lambda}(r) = \lambda\psi(r/\lambda), \quad r\ge \lambda.
$$%
Let $u$ be a $C^{2}(\Omega)\cap C^{0}(\overline\Omega)$ solution of~\ref{sme}, and let 
$$%
\beta=\sup_{r>1}\psi(r)/r \text{ (which is the same as $\beta=\sup_{r>\lambda}\psi_{\lambda}(r)/r$ for each $\lambda>0$)}.
$$%
We remark that actually $\beta=(\frac{m-1}{n-1})^{1/2}$ if $m+n-1\ge 7$, because the solution $\psi(r)$ remains below
$(\frac{m-1}{n-1})^{1/2}r$ for all $r$ in this case; but in any case for any $m,n\ge 2$ such a finite $\beta$ exists because
$\frac{\psi(r)}{r}\to (\frac{m-1}{n-1})^{1/2}$ as $r\uparrow\infty$.  Now assume $u=\varphi$ on $\partial\Omega$, where
\ab{3pt}{3pt}{%
$$%
\varphi(x) >\sup\{\beta|x-x_{0}|:x_{0}\in\Omega\}\text{ for each $x$  }\in \partial\Omega. %
\dl{st-pos-ex}
$$}%
Then evidently $\varphi(x)\ge \psi_{\lambda}(|x-x_{0}|)$ for each $(x,\lambda)\in \partial\Omega\times(0,\infty)$ with
$|x-x_{0}|>\lambda$, and we can select the largest $\lambda>0$ such that $u(x)\ge\psi_{\lambda}(|x-x_{0}|)$ for every $x\in
\Omega\cap \{y:|y-x_{0}|>\lambda\}$, and, for such a $\lambda$, we have $\xi\in \Omega\cap \{x:|x-x_{0}|>\lambda\}$ with
$u(\xi)=\psi_{\lambda}(|\xi-x_{0}|)$ and $u(x)\ge \psi_{\lambda}(|x-x_{0}|)$ for all $x$ in some neighborhood of
$\xi$.  But then, by taking the difference of the SME for $u(x)$ and the SME for $\psi_{\lambda}(|x-x_{0}|)$, we see that
then $v(x)=u(x)-\psi_{\lambda}(|x-x_{0}|)$ has a local minimum value of zero at $x=\xi$ and in a neighborhood of $\xi$
it satisfies an equation of the form $\sum_{i,j=1}^{n}a_{ij}D_{i}D_{j}v+\sum_{i=1}^{n}b_{i}D_{i}u+cu=0$, with
$(a_{ij})$ positive definite and $a_{ij},b_{i}, c$ continuous. This evidently contradicts the Hopf maximum principle. So
we have proved that any $C^{2}(\Omega)\cap C^{0}(\overline\Omega)$ solution of \ref{sme} with $u(x)\ge
\sup\{\beta |x-x_{0}|:x_{0}\in\Omega\}$ for each $x\in \partial\Omega$
automatically satisfies $u(x)\ge \sup\{\psi_{\lambda}(|x-x_{0}|):x_{0}\in\Omega,\lambda\in (0,\infty)\}$  for each $x\in
\Omega$.  Since $\sup_{\lambda\in (0,\infty)}\psi_{\lambda}(|x-x_{0}|)\ge (\frac{m-1}{n-1})^{1/2}|x-x_{0}|$ we see
that  
$$%
\inf_{x\in\partial\Omega}\sup_{x_{0}\in\Omega,\lambda\in (0,1)}\psi_{\lambda}(|x-x_{0}|)
\ge  \inf_{x\in\partial\Omega}\sup_{x_{0}\in\Omega}\bigl(\tfrac{m-1}{n-1}\bigr)^{1/2}|x-x_{0}| \ge \gamma 
$$%
for suitable $\gamma=\gamma(m,\Omega)>0$, so indeed any continuous data $\varphi$ satisfying \ref{st-pos-ex} is
strongly positive.

\smallskip

Now, we can state the main existence result of this section

\begin{state}{\tl{existence2}\,Theorem.}%
  Let $\Omega$ be a uniformly convex $C^{2,\alpha}$ domain in $\R^{n}$, $n\ge 2$, and
$\{\varphi_{\lambda}\}_{\lambda\in [0,1]}\subset C^{2,\alpha}(\overline\Omega)$ such that
$(x,\lambda)\to\varphi_{\lambda}(x)$ is $C^{0}$ map $\overline\Omega\times [0,1]\to\R$, $\varphi_{1}>\varphi$, where
$\varphi$ is any strongly positive data (as in~\emph{\ref{st-pos-non-solv}}), and $\varphi_{0}$ non-solvable (also as
in~\emph{\ref{st-pos-non-solv}}).  Then there is $\lambda\in (0,1)$ such that there is a singular solution $u$ of
\emph{\ref{sme}} with $u=\varphi_{\lambda}$ on $\partial \Omega$ and $\dist(\sing u,\partial\Omega)>0$.
\end{state}

{\bf{}Remark:} With some more effort it is possible to replace the uniform convexity hypothesis with the hypothesis that
$\partial\Omega$ is  mean convex, but we shall not discuss that here.

\smallskip

The proof of Theorem\thn\ref{existence2} involves an application of a standard Leray-Schauder degree argument, aided
by the H\"{o}lder continuity results of the previous section.  The Leray-Schauder component of the proof is presented in
the following lemma:

\smallskip

\begin{state}{\tl{schauder}\,Lemma.}%
  Let $V$ be an  open (not necessarily bounded)  subset of a Banach space $\mathcal{B}$ 
and let $T_{\lambda}: \overline V\to {\cal{}B}$,
$0\le \lambda\le 1$ be such that the map $(x,\lambda)\mapsto T_{\lambda}(x), (x,\lambda)\in
\overline{V}\times [0,1]$, is a continuous compact map (i.e.\ a continuous map taking bounded subsets of
$\overline V\times[0,1]$ into  compact subsets of ${\cal{}B}$), and assume %

\emph{(a)\,\,} $T_{1}$ is a constant map $\overline V\to \mathcal{B}$ with constant value $p_{0}\in V$ \hfill\break %
\emph{(b)\,\,} $T_{0}$ has no fixed points in $\overline V$ \hfill\break %
\emph{(c)\,\,} $\sup\{ \|u\|: u\in \cup_{\lambda\in [0,1]}\{v\in \overline{V}:  T_{\lambda}(v)=v \}\}<\infty$.

Then there is a  $u \in \partial V$ and a $\lambda\in (0,1)$ with $T_{\lambda}(u)=u$.
\end{state}

\begin{proof}{\bf{}Proof of Lemma \ref{schauder}:} The proof is a standard application of the Leray-Schauder degree of
completely continuous maps (i.e.\ maps of the form $\iota -T$, where ${T}$ is continuous and compact and where
$\iota$ is the identity map on ${\cal{}B}$).   Specifically we use the fact (see e.g.\ \cite{Dei85})  that 
 if $U$ is a bounded open subset of a Banach space $\mathcal{B}$, then there is a well-defined
topological degree $d$ for completely continuous transformations of $\overline U$ into $\mathcal{B}$ as follows:

\smallskip

(i) If $T: \overline U\to \mathcal{B}$ is continuous and compact then the topological degree $d(\iota -T,U,q)$ is a
well-defined integer for $q\in \mathcal{B}\setminus (\iota -T)(\partial U)$, $d(\iota -T,U,q)$ remains constant for $q$ in
a given connected component of $\mathcal{B}\setminus (\iota -T)(\partial U)$, and $d(\iota -T,U,q)\neq 0\Rightarrow 
q\in (\iota-T)U$.

\smallskip

(ii) If $T_{\lambda}: \overline U\to \mathcal{B}$ are given for $\lambda \in [0,1]$ such that the map
$(p,\lambda)\mapsto T_{\lambda}(p),\,(p,\lambda) \in \overline U\times [0,1]$, is a compact continuous map, and if
$q\in \mathcal{B}\setminus (\cup_{\lambda\in [0,1]}(\iota -T_{\lambda})(\partial U))$, then $d(\iota -T_{\lambda}, U,
q)= d(\iota -T_{0}, U, q)$ for each $\lambda\in [0,1]$.

\smallskip

(iii) If $T$ is a constant map with constant value $q_{0}\in U$, then $d(\iota -T,U,0)=1$.

\smallskip

To prove the lemma, we first use hypothesis~(c) to choose $R>0$ such that
$$%
R> \sup_{\{u: u\in   \overline{V} \text{ and } T_{\lambda}(u)=u  %
                                                                    \text{ for some }\lambda\in   [0,1]\}} \|u\|,
$$%
and then we apply the above properties~(i), (ii), (iii) of the topological degree with $U=V\cap \{u\in \mathcal{B}: \|u\|<R\}$
as follows:

Either $\exists\, \lambda\in [0,1]$ and $u\in \partial (V\cap \{u:\|u\|<R\})(\subset(\partial{V})\cup\{u:\|u\|=R\})$ with
$T_{\lambda}(u)=u$ or else there is no such $\lambda$. But by property~(ii) of the degree the latter alternative implies
$d(\iota -T_{\lambda},V\cap \{u: \|u\|<R\}, 0)$ is constant for $\lambda\in [0,1]$. By hypothesis~(a) and property~(iii) the
constant value must be~1, and hence, by property~(i), $T_{0}(u)=u$ for some $u\in V\cap \{u:\|u\|<R\}$, contradicting
hypothesis~(b) of the lemma.  Thus the former alternative holds, and, since $\|u\|<R$ whenever $T_{\lambda}(u)=u$ with
$u\in \overline{V}$ and $\lambda\in [0,1]$, we deduce that there is a $\lambda$ with $T_{\lambda}(u)=u$ for some $u\in
\partial V$. \end{proof}

\begin{proof}{\bf{}Proof of Theorem\thn\ref{existence2}:}
  Let $\delta\in (0,1]$.  We apply the Lemma\thn\ref{schauder} with ${\cal{}B}=C^{1,\alpha}(\overline\Omega)$ and $V
=\{u \in C^{1,\alpha} (\overline{\Omega}) : u >\delta\}$. Let $a_{ij}(p)=\delta_{ij}-p_ip_j/(1+|p|^2)$, so that
\ab{3pt}{3pt}{%
$$%
\mathcal{M}(u)=\frac{1}{\sqrt{1+|Du|^2}}\,a_{ij}(Du)D_{ij}u,
$$}%
and consider the following family of problems for given $v \in \mathcal{S}_\delta$:
$$%
Q_\lambda^v=\left\{
\begin{aligned}
      &a_{ij}(Dv)D_{ij}u=\frac{(m-1)}{v},  %
     \,\,u_{|\partial \Omega}=\varphi_{2\lambda}, & 0 \leq \lambda \leq \frac{1}{2} \\ 
    &a_{ij}(Dv)D_{ij}u=\frac{(m-1)}{v},  %
     \,\,u_{|\partial \Omega}=\varphi_1+4(\lambda -\frac{1}{2})(K-\varphi_1),  %
                                               & \frac{1}{2} < \lambda \leq \frac{3}{4}\\     
    &a_{ij}(Dv)D_{ij}u=\frac{4(m-1)(1-\lambda)}{v},\,\,  u_{|\partial \Omega}=K, & \frac{3}{4} < \lambda \leq 1
\end{aligned}   \right. 
$$%
where $K$ is to be chosen (large).

Each of these equations is a  inhomogeneous linear second order elliptic equation with H\"older continuous
coefficients and with no first or zero order terms, so by the  theory of such equations (see \cite{GT})
each problem has a unique solution in $C^{2,\alpha}(\overline{\Omega})$.  Define the map
$T_\lambda : \overline{V} \mapsto C^{1,\alpha}(\overline{\Omega})$ by
\ab{3pt}{3pt}{%
$$%
T_\lambda (v)=u,
$$}%
where $u$ is the solution to the problems above.  As the coefficients of $Q^v_\lambda$ are in
$C^{0,\alpha}(\overline\Omega)$, by virtue of global Schauder estimates in \cite[Chapter 6]{GT} $T_\lambda$ maps
bounded sets in $\overline V$ to bounded sets in $C^{2, \alpha}(\overline{\Omega})$, which are precompact in
$C^{1,\alpha}(\overline{\Omega})$.  Therefore, $T_\lambda$ is a compact mapping. The continuity and compactness of
$T:\overline{V}\times[0,1]\to C^{1,\alpha}(\overline\Omega)$ (with $T(u,\lambda)=T_{\lambda}(u)$) follows from a
similar standard argument.  Also by the maximum principle for elliptic equations, $T_1\equiv K$. By assumption
$T_0$ has no fixed points in $\overline{V}$.  Using exactly the same barrier argument discussed in [GT,
Chapter\thn14] (capitalizing on the fact that the fixed points satisfy ${\cal{}M}u=\frac{m-1}{u\sqrt{1+|Du|^{2}}}$ and
$0<\frac{m-1}{u\sqrt{1+|Du|^{2}}}\le \delta^{-1}\frac{m-1}{\sqrt{1+|Du|^{2}}}$ for $u\in \overline{V}$) the boundary
gradient estimate
$$%
\sup_{\partial \Omega} |Du| \leq C, %
\pdl{bar1}
$$%
holds for the fixed points of $T_\lambda$, where $C=C(m,n,\delta,  \sup_{\mu}|\varphi_{\mu}|_2)$.

Then by applying the gradient estimates of~\cite{Sim76} in the case when the functions $A_{j},B$ of \cite{Sim76} satisfy
\[%
A_{j}(x,u,Du)=D_{j}u/\sqrt{1+|Du|^{2}}\text{ and }|B(x,u,Du)|\le C/\sqrt{1+|Du|^{2}}, 
\]%
we have
\ab{3pt}{3pt}{%
$$%
\sup_{\Omega}|Du| \le C, \quad C=C(M,\Omega,\delta),
\pdl{bar3}
$$}%
where  $M$ is any upper bound for $\sup_{\lambda\in [0,1]}|\varphi_{\lambda}|_{C^{2}(\overline\Omega)}$.
(Alternatively one can use the argument of \cite{DieH90} for this.)

Combining  \ref{bar3} and the regularity theory for quasilinear elliptic equations in \cite{GT} ensures that
$$%
  \sup_{\{u \in \overline{V} : T_\lambda (u)=u \text{ for some } \lambda \in [0,1]\} } \|u\| <
  \infty. %
  \pdl{bar4}
$$%
Now we can apply Lemma \ref{schauder} to conclude that there is $u_{\delta} \in \partial V$ and $\lambda_{\delta} \in
(0,1]$ with $T_{\lambda_{\delta}}(u_{\delta})=u_{\delta}$.

\medskip

Next we claim that for suitably large $K$ (depending on $\delta$) we can arrange that
\ab{3pt}{3pt}{%
$$%
\begin{aligned}%
&\text{If $\delta<\eta$ with $\eta=\eta(\varphi_{1})>0$ as in~\ref{st-pos-non-solv} then}  \\  %
\noalign{\vskip-3pt}
&\text{\hskip1in there are no fixed points of $T_\lambda$ on $\partial   V$ for $\lambda \in
[\tfrac{1}{2}, 1]$} %
 \end{aligned}%
\pdl{bar5}
$$}%
(so that all fixed points of $T_{\lambda}$ occur for $\lambda\in (0,\frac{1}{2}]$).

To see~\ref{bar5},  first consider the case $\lambda \in [\frac{1}{2}, \frac{3}{4}]$. Assuming $K > \sup \varphi_1$, we have
$$%
\varphi_1 + 4 (\lambda -\tfrac{1}{2})(K-\varphi_1) \geq \ \varphi_{1},
$$%
and so $u\ge \eta>\delta$ by \ref{st-pos-non-solv}, contradicting the fact that $u\in\partial V$.

Next we consider the case where $\lambda \in (\frac{3}{4}, 1)$. Recall that $a_{ij}(p)=\delta_{ij}-p_ip_j/(1+|p|^2)$, so we
have
$$%
\tsum_{i=1}^{n} a_{ii} > n-1.
\pdl{bar6.5}
$$%
Supposing that there is a fixed point $u\in\partial V$,  after dividing by $K$ in \ref{sme} we get
$$%
\tsum a_{ij}(Du)D_{ij}\bigl(\frac{u}{K}\bigr)=\frac{4(1-\lambda)(m-1)}{uK}\;,\;u \in \partial V.
\pdl{bar6.6}
$$%
Assume without loss of generality that $0\in\Omega$ and let $\theta = \frac{1}{2}\mbox{diam}^{-2}(\Omega)$. Then
$w(x)=\frac{u}{K} -\theta |x|^2$ satisfies
$$%
w_{|\partial \Omega} \geq \tfrac{1}{2}  \quad \text{(since $u\ge K$ on $\partial\Omega$)}
$$%
and, by~\ref{bar6.5} and \ref{bar6.6},
$$%
\tsum a_{ij}(Du)D_{ij}w < \tfrac{4(1-\lambda)(m-1)}{uK}-2(n-1)\theta. 
$$%
Thus if $K > \max \Big\{\frac{2(n-1)\theta}{(m-1)\delta}, \sup \varphi_1\Big\}$ then $a_{ij}(Du)D_{ij}w < 0 \text{ in
}\Omega$, so the maximum principle implies that
$$%
\inf_{\Omega} w\geq \tfrac{1}{2}.
$$%
In particular, $\inf_{\Omega} u \ge K/2>\delta$ contradicting the fact that $\inf_{\Omega}u=\delta$.  Thus \ref{bar5} is
established.

On the other hand Lemma~\ref{schauder} implies there is indeed a fixed point in $\partial V$ of
$T_{\lambda}$ for some $\lambda\in (0,1)$, and \ref{bar5} ensures that this $\lambda$ is in the interval $(0,1/2]$, which
means that $u$ satisfies~\ref{sme} with $u|\partial\Omega=\varphi_{\mu}$ for some $\mu\in (0,1]$.  Taking
$\delta=\delta_{j}<\inf_{\mu\in [0,1]}\min_{\partial\Omega}\varphi_{\mu}$ with $\delta_{j}\to 0$, we thus see that
there are $u_{j}$ such that $u_{j}$ satisfies~\ref{sme} and $u_{j}=\varphi_{\mu_{j}}$ for
some $\mu_{j}\in (0,1)$ and $\min_{\Omega}u_{j}=\delta_{j}$.

We claim that there is $\eta>0$ such that
$$%
u_{j}(y)\ge \eta \text{ for all $j$ and all $y\in\{x\in\Omega:\dist(x,\partial\Omega)<\eta\}$.} %
\pdl{tag8}
$$%
This is easily checked using the natural variant of the argument, involving the Hopf maximum principle and the family
$\psi_{\lambda}$ of ODE solutions of~\ref{sme}, mentioned earlier in this section, as follows:

Let $x_{0}\in \partial\Omega$, $\nu_{0}$ be the inward pointing unit normal of $\partial\Omega$ at $x_{0}$, and take a
point $y_{0}=x_{0}+t_{0}\eta_{0}$ on the ray $\{x_{0}+t\eta_{0}:t>0\}$ with $t_{0}$ sufficiently large to ensure that the
sphere $S_{y_{0}}$ with center $y_{0}$ and radius $t_{0}$ satisfies $S_{y_{0}}\cap \partial\Omega=\{x_{0}\}$.  Then
$\psi_{\lambda}(|x-y_{0}|)$ is a family of solutions defined over the region $|x-y_{0}|\ge \lambda$ and, if $\epsilon>0$ is
small enough, each $\lambda\in(t_{0}-\epsilon,t_{0})$ gives a solution of~\ref{sme} which satisfy $\psi_{\lambda}\le
\inf_{\mu\in [0,1]}\varphi_{\mu}$ on $\{x\in
\partial\Omega:|x-y_{0}|\ge \lambda\}$ for each $\lambda$ in the interval $(t_{0}-\epsilon,t_{0})$.  We then claim that the
solutions $u_{j}$ constructed above must satisfy, for each $\lambda$ in the interval $(t_{0}-\epsilon,t_{0})$,
$$%
u_{j}(x)\ge\psi_{\lambda}(|x-y_{0}|), \quad \forall\,x\text{ with }  |x-y_{0}|\ge t_{0}-\lambda.  %
\pdl{tag9}
$$%
Otherwise we pick $\lambda_{j}$ to be the smallest $\lambda$ in the interval $(t_{0}-\epsilon,t_{0})$ where this is true. 
Then, analogous to the argument used above to prove the existence of strongly positive data, we contradict the Hopf
maximum principle for the difference $u_{j}(x)-\psi_{\lambda_{j}}(|x-y_{0}|)$.  So, since~\ref{tag9} applies for each
point $x_{0}\in\partial\Omega$,  \ref{tag8} is proved.

Now in the region $\{x\in\Omega:\dist(x,\partial\Omega)\ge\eta\}$ we have equicontinuity of the $u_{j}$ by
Corollary~\ref{holdercor}, and on the other hand in the region $\{x\in\Omega:\dist(x,\partial\Omega)<\eta\}$ we obtain
the same local gradient estimates (using \cite{Sim76}) as mentioned above (since $u_{j}\ge \eta$ in this region), hence we
also obtain equicontinuity of the $u_{j}$ $\{x\in\Omega:\dist(x,\partial\Omega)<\eta\}$.  So there is a subsequence
$u_{j'}$ converging uniformly on $\overline\Omega$ to a singular solution $u$ of~\ref{sme} with $\sing
u\cap\{x\in\Omega:\dist(x,\partial\Omega)\le\eta\}=\emptyset$ and $u=\varphi_{\lambda}$ for some $\lambda\in (0,1)$. 
This completes the proof of Theorem~\ref{existence2}.
\end{proof}

\section{Compactness and Regularity Results}\label{sec-regularity}

In this section we discuss the regularity of limits of the graphs of the solutions of the SME (i.e.\ equation~\ref{sme}) in
the open upper half space $x_{n+1} >0$ by showing that the regularity theory of~\cite{SS81} can be applied.

We let
\begin{align*}%
&{\cal{}G}_{R}=\{G(u): \text{ with $u>0$ a $C^{2}$ solution of \ref{sme} on $B_{R}^{n}$}\},\\
\noalign{\vskip-3pt and}
&\overline{{\cal{}G}}_{R} = \{ \text{varifolds $V$ on $B_{R}^{n}\times \R$ expressible as } %
V=\lim G_{j} \text{ in }B_{R}^{n}\times\R  \\
\noalign{\vskip-2pt}
&\hskip1.5in \text{for some sequence $u_{j}\in {\cal{}G_{R}}$, where $G_{j}=G(u_{j})\in {\cal{}G}_{R}$}\}.
\end{align*}%
Of course here $G_{j}$ is interpreted as the multiplicity~$1$ varifold corresponding to the graph $G(u_{j})$ of $u_{j}$
and the varifold convergence guarantees the measure convergence of ${\cal{}H}^{n}\res\, G_{j}$ to a Radon measure on
$B^{n}_{R}\times\R$. In view of the volume bounds of~\ref{vol-bds}(i), the compactness theorem for Radon measures on
the space $G_{n}(B^{n}_{R}\times\R)=(B_{R}^{n}\times\R)\times G(n,n+1)$ (where $G(n,n+1)$ denotes the
$n$-dimensional subspaces in $\R^{n+1}$ with the metric $\rho(S,T)=|p_{S}-p_{T}|$, where $p_{S},p_{T}$ the orthogonal
projections of $\R^{n+1}$ onto the $n$-dimensional subspaces $S,T$) guarantees that every sequence $G_{j}=G(u_{j})$
automatically has a subsequence $G_{j'}$ with $\lim G_{j'}=V$ for some varifold $V$---i.e.\ for some Radon measure
$V$ on $G_{n}(B^{n}_{R}\times\R)$. For further discussion of the theory of varifolds we refer e.g.\ to \cite[Chapter
8]{Sim83}.

We claim that 
$$%
{\cal{}H}^{n}(\spt V\cap (B^{n}_{R}\times\{0\})) = 0  \text{ for each $V\in\overline{\cal{}G}_{R}$}. %
\dl{sing-finite-meas}
$$%
To prove this, suppose $V=\lim G(u_{j})$ with $u_{j}>0$ a solution of~\ref{sme} on $B_{R}^{n}$, and
on the contrary that ${\cal{}H}^{n}(K)>0$, where
$$%
K=\spt V\cap (B^{n}_{R}\times\{0\}).
$$%
We claim $u_{j}(x)\to 0$ for ${\cal{}H}^{n}$-a.e.\ $(x,0)\in K$. Otherwise $\limsup u_{j}(x)>0$ on a subset of $K$ with
positive ${\cal{}H}^{n}$ measure, hence there would be $\epsilon_{0}>0$ and a set $K_{0}\subset K$ with
${\cal{}H}^{n}(K_{0})>0$ and $\limsup u_{j}(x)\ge \epsilon_{0}$ for every $x\in K_{0}$. Also for every $\sigma>0$ we
must have some $\xi_{j}\in G_{j}\cap (B_{\sigma}^{n}(x)\times[0,\epsilon_{0}/2])$ for all sufficiently large $j$ (otherwise,
since $G_{j} \to V$ in the varifold sense on $B_{R}\times\R$, we would have $(x,0)\notin\spt V$).  So then
$\{u_{j}(tx+(1-t)\xi_{j}):t\in [0,1]\}$ includes every value between $\epsilon_{0}/2$ and $3\epsilon_{0}/4$ for infinitely
many $j$.  Thus for each $\sigma>0$ there are infinitely many $j$ with $G_{j}\cap
(B^{n}_{\sigma}(x)\times\{t\})\neq\emptyset$ for every $t\in [\epsilon_{0}/2,3\epsilon_{0}/4]$, so in particular, since
$G_{j}\to V$ locally in the Hausdorff distance sense on $B_{R}^{n}\times(0,\infty)$, the vertical segment $\{(x,t): t\in
[\epsilon_{0}/2,3\epsilon_{0}/4]\}$ is contained in $\spt V$. So $K_{0}\times[\epsilon_{0}/2,3\epsilon_{0}/4]\subset \spt
V$ contradicting the fact that $\spt V\cap(B^{n}_{\rho}\times(0,\rho))$ has finite $n$-dimensional measure for each
$\rho<R$ by the volume bounds of Lemma~\ref{vol-bds}(i).

So indeed $u_{j}(x)\to 0$ for ${\cal{}H}^{n}$-a.e.\ $(x,0)\in K$ and by Egoroff's theorem there is then a subset $S\subset
K$ of positive measure with $u_{j}\to 0$ uniformly on $S$.  Pick $(y,0)\in S$ and $\rho>0$ with
$B^{n+1}_{\rho}(y,0)\subset B^{n}_{R}\times\R$ and ${\cal{}H}^{n}(S\cap B^{n+1}_{\rho/4}(y,0))>0$. Then, with
$\widetilde S=S\cap B_{\rho/4}(y,0)$, by Cauchy-Schwarz and the fact that $\{(x,u_{j}(x)):x\in \widetilde S\}\subset
B^{n+1}_{\rho/2}(y,u_{j}(y))$ for sufficiently large $j$,
\begin{align*}%
  \bigl({\cal{}H}^{n}(\widetilde S)\bigr)^{2}  %
  &= \bigl(\tint_{\widetilde S}   (1+|Du_{j}|^{2})^{1/4}(1+|Du_{j}|^{2})^{-1/4}\bigr)^{2}    \\  %
  &\le \tint_{\widetilde   S}(1+|Du_{j}|^{2})^{1/2}  \tint_{\widetilde S}(1+|Du_{j}|^{2})^{-1/2}  \\  %
  & \le C\tint_{\widetilde S}(1+|Du_{j}|^{2})^{-1/2} \qquad \text{(by \ref{vol-bds}(i))} \\  %
  & \le C{\sup}_{\widetilde S\,}u_{j}\, \tint_{\widetilde S}u_{j}^{-1}(1+|Du_{j}|^{2})^{-1/2} \\ %
  & \le C{\sup}_{S}u_{j}\, \tint_{B^{n}_{\rho/2}(y)}u_{j}^{-1}(1+|Du_{j}|^{2})^{-1/2}  %
\end{align*}%
for all sufficiently large $ j$.  On the other hand by using the weak version \ref{sme}$^{\prime}$ of~\ref{sme} with
$u_{j}$ in place of $u$ and with $\zeta$ a cut-off function with $\zeta= 1$ on $B_{\rho/2}(y)$ and $\zeta=0$ outside
$B_{\rho}(y)$, we conclude that
$$%
\tint_{B^{n}_{\rho/2}(y)}u_{j}^{-1}(1+|Du_{j}|^{2})^{-1/2} \le C, %
$$%
with $C$ independent of $j$. Combining the previous two inequalities we have
$$%
({\cal{}H}^{n}(\widetilde S))^{2}\le    C\hskip0.5pt {\sup}_{S}u_{j}  \to 0,
$$%
a contradiction. So~\ref{sing-finite-meas} is proved.

We also define
$$%
\overline{\cal{}G}_{\infty} =\text{ the set of varifolds $V$ on $\R^{n}\times\R$ with } %
V\res\, (B^{n}_{R}\times\R)\in\overline{\cal{}G}_{R} \,\forall\,R>0.
$$%
As already mentioned above, by the volume bounds of Lemma~\ref{vol-bds} we can use the compactness theorem for
Radon measures to prove that an arbitrary sequence $G(u_{j})\in{\cal{}G}_{R}$ always has a subsequence which converges
to $V\in \overline{\cal{}G}_{R}$ in the sense of Radon measures in $(B^{n}_{R}\times\R)\times G(n,n+1)$.  Additionally,
since the mean curvature of $G_j$ is bounded by $C/\sigma$ in the region $\{x_{n+1} >\sigma\}$ ($\sigma >0$
arbitrary), by the Allard compactness theorem any $V\in\overline{{\cal{}G}}_{R}$ is an integer multiplicity varifold with
(generalized) mean curvature bounded by $C/\sigma$ in $B_R^{n}\times (\sigma,\infty)$ for each $\sigma>0$.

\smallskip

We proceed to show that the regularity theory of~\cite{SS81} can be applied locally in $B^{n}_{R}\times(0,\infty)$. First
take an arbitrary $C^{3}(\Omega)$ function $u$ (not necessarily satisfying~\ref{sme}), and let $G=G(u)$, $H=$ mean
curvature of $G=\sum_{i=1}^{n}D_{i}(D_{i}u/\sqrt{1+|Du|^{2}})$,  $w=-\log \nu_{n+1}$,
$\nu_{n+1}=1/\sqrt{1+|Du|^{2}}$ (so that $w\ge 0$ because $0<\nu_{n+1}\le 1$). \ref{laplace} in the proof of
Theorem\thn\ref{holder1} says 
$$%
\Delta_{G}\nu_{n+1} + |A_{G}|^{2}\nu_{n+1} = -e_{n+1}\cdot \nabla_{G}H,
\dl{limits0} 
$$%
and (using $\nabla_{G}w=-\nu_{n+1}^{-1}\nabla_{G}\nu_{n+1}$) the correspond identity for $\Delta_{G}w$:
$$%
\Delta_{G}w-(|A_{G}|^{2}+|\nabla_{G}w|^{2})=\nu_{n+1}^{-1}e_{n+1}\cdot \nabla_{G}H,
\dl{limits1}
$$%
where $A_{G}$ is the second fundamental form of $G$ and $\nabla_{G}$ denotes the gradient on $G$, so that
$$%
|A_{G}|^{2} = \nu_{n+1}^{2}\tsum_{i,j=1}^{n}g^{i\ell} g^{jk}D_{i}D_{j}u D_{k}D_{\ell}u,
$$%
where $g^{ij}=\delta_{ij}-\nu_{i}\nu_{j}$.

When $u$ is a regular or singular solution of \ref{sme} we have, in $G(u)\cap (B_{R}^{n}\times(0,\infty))$,
$H=(m-1)\nu_{n+1}/u$ (which can be expressed $H=(m-1)\nu_{n+1}/x^{n+1}$), and so the identities~\ref{limits0}
and~\ref{limits1} can be written
\begin{align*}%
\tg{limits1.5}
  &\Delta_{G}\nu_{n+1} + |A_{G}|^{2}\nu_{n+1} = -(m-1)e_{n+1}\cdot \nabla_{G}((x^{n+1})^{-1}\nu_{n+1}) \\
  \tg{limits2}  &\Delta_{G}w-(|A_{G}|^{2}+|\nabla_{G}w|2)=  \\
&\hskip0.5in    -(m-1)(x^{n+1})^{-2}e_{n+1}\cdot \nabla_{G}x^{n+1} -(m-1)(x^{n+1})^{-1}e_{n+1}\cdot \nabla_{G}w.  %
\end{align*}%
Rewriting \ref{limits2} in the weak form, we have
\begin{align*}%
\tg{limits3}&  \tint_{G} (|\nabla_{G} w|^2+|A|^2)\zeta^2  = -\tint_{G} 2\zeta \nabla_{G}w.\nabla_{G} \zeta   \\
\noalign{\vskip-1pt} %
&\hskip0.4in
  -(m-1)\tint_{G}\bigl( (x^{n+1})^{-2}e_{n+1}\cdot \nabla_{G}x^{n+1} %
                         +(x^{n+1})^{-1}e_{n+1}\cdot  \nabla_{G}w\bigr)\zeta^{2}
\end{align*}%
provided $\zeta\in C^{1}_{c}(\R^{n+1})$ with $\spt \zeta \cap G$ is a compact subset of $G$.  Using the Cauchy-Schwarz
inequality we have $2\zeta \nabla w.\nabla_{G} \zeta\le (1-\epsilon) \zeta^{2}|\nabla_{G}w|^{2}
+(1-\epsilon)^{-1}|\nabla_{G}\zeta|^{2}$ and $(x^{n+1})^{-1}e_{n+1}\cdot \nabla_{G}w\bigr)\le
\epsilon|\nabla_{G}w|^{2}+\epsilon^{-1}(x^{n+1})^{-2}$, so \ref{limits3} gives
$$%
\tint_{G} |A|^2\zeta^2 \leq %
\frac{1}{1-\epsilon}\tint_{G} |\nabla_{G} \zeta|^2 + \epsilon^{-1}\delta^{-2}\tint_{G}\zeta^2 %
\dl{ssf}
$$%
provided $\spt\zeta\cap G$ is compact and contained in the region $x^{n+1}\ge \delta$.

The inequality~\ref{ssf} enables us to use the regularity theory of \cite{SS81} for solutions $u$ in the region $x^{n+1}>0$,
which we now discuss:

Define $F: \mathbb{R}^{n+1}_+ \times \mathbb{R}^n \rightarrow \mathbb{R}$ by
$$%
F(x,x_{n+1},p)=x_{n+1}^{m-1}|p|
$$%
With this notation the Symmetric Minimal Surface equation \ref{sme} can be viewed as the Euler- Lagrange equation for
the functional
$$%
\mathcal{F}(G)=\tint_{G} F(x,x^{n+1}, \nu)\,d\mathcal{H}^n
$$%
where $\nu=(-Du,1)/\sqrt{1+|Du|^{2}}$ is the upward pointing unit normal  on $G=G(u)$.  If we restrict the domain to
where $x_{n+1} >\delta$, for some fixed $\delta >0$, then $F$ satisfies all the properties (1.2) to (1.6) in \cite{SS81}. 
Note that, in the region $x_{n+1}>0$, the inequality \ref{ssf} replaces the $F$-stability in \cite{SS81}, because \ref{ssf}
implies the stability inequality \cite[(1.17)]{SS81} required in the proof of regularity results in \cite{SS81}; in fact \ref{ssf}
has an additional factor $(1-\epsilon)^{-1}$ in front of the principal term on the right side, but this causes no
complication in the discussion of \cite{SS81} so long as we take $\epsilon=\epsilon(n)$ sufficiently small.

Thus with $V=\lim G_{j}\in\overline{\cal{}G}_{R}$ as above, we can apply the regularity and compactness theory
of~\cite{SS81} and in particular Theorem 4 in \cite{SS81} implies that $V$ has a singular set in $x_{n+1}>0$ (i.e.\ $\sing
V\cap(B_{R}\times(0,\infty))$) of codimension at least 7 (empty for $n\le 6$, discrete for $n=7$).

We henceforth write
$$%
\Sigma = \reg V\cap (B^{n}_{R}\times(0,\infty)),
$$%
and observe that, in view of the regularity statements above,
$$%
\dim \sing \Sigma \le \max\{n-7,0\},
$$%
where $\sing \Sigma\cap (B^{n}_{R}\times(0,\infty)) = (\overline \Sigma\setminus\Sigma) \cap
B_{R}^{n}\times(0,\infty)$.

\smallskip

An additional part of the theory established in \cite{SS81} shows that, in the region where $x^{n+1}>0$, the
approximating graphs $G_{j}$ actually converge in the $C^{k}$ sense to $\reg V$ for each $k$ locally near points of
$\reg V$. More precisely:

 \begin{state}{\tl{sheeting}\,Lemma (``Sheeting Lemma.'')}%
   Suppose that $G_j$ is the sequence of graphs of (possibly singular) solutions $u_{j}$ of~\emph{\ref{sme}} converging to
$V$ is the sense of varifolds in the cylinder $B_{R}^{n}\times(0,\infty)$, let $\Sigma=\reg V\cap
(B_{R}^{n}\times(0,\infty))$ and let $(\xi,\tau)\in \Sigma$.  Then there is $\rho>0$ and an integer $j_{0}$ such that for all
$j\ge j_{0}$ there is an integer $L\ge 1$ and $C^{\infty}$ functions $u^1_j <\cdots <u^L_j$ on $\Sigma\cap
B^{n+1}_{\rho}(\xi,\tau)$ such that
$$%
G_{j} \cap B^{n+1}_{\rho/2}(\xi,\tau) =\cup_{i=1}^L G_{\Sigma}(u_i^j) \cap B^{n+1}_{\rho/2}(\xi,\tau),
$$%
and $\sup_i |u_i^j|_{C^k(\Sigma\cap B^{n+1}_{\rho}(\xi,\tau))}\rightarrow 0$ as $j\rightarrow \infty$ for each
$k=0,1,2,\ldots$. Here $G_{\Sigma}(u_{j}^{i})$ denotes the graph of $u_{j}^{i}$ defined by
$$%
G_{\Sigma}(u_{j}^{i}) = \{x+u_{j}^{i}(x)\nu(x):x\in \Sigma\cap B^{n+1}_{\rho}(\xi,\tau))\},
$$%
where $\nu$ is a smooth unit normal for $\Sigma$.
\end{state}

\medskip

Of course in view of the $C^{k}$ convergence of the above lemma, we can take limits in \ref{ssf} in order to deduce
the inequality
$$%
\tint_{\Sigma} |A_{\Sigma}|^2\zeta^2 \leq %
\frac{1}{1-\epsilon}\tint_{\Sigma} |\nabla_{\Sigma} \zeta|^2 + \epsilon^{-1}\delta^{-2}\tint_{\Sigma}\zeta^2 %
\dl{ssf1}
$$%
for each $\epsilon,\delta>0$ and each $\zeta\in C^{1}_{c}(B_{R}^{n}\times(\delta,\infty))$, where  $A_{\Sigma}$
denotes the second fundamental form of $\Sigma$ and $\nabla_{\Sigma}$ is the gradient operator on $\Sigma$.

As a further consequence of Lemma~\ref{sheeting} we can prove that $V$ has multiplicity $N\le 2$ at all points of
$\Sigma$:

 \begin{state}{\tl{mult2}\,Lemma.}%
   The varifold $V$ has multiplicity $\le 2$ at each point of $\Sigma=\reg V\cap (B_{R}^{n}\times(0,\infty))$.
\end{state}

{\bf{}Remark:}\label{const-mult} Notice that, since $\reg V$ is an embedded $C^{2}$ submanifold, by the constancy
theorem for integer multiplicity varifolds with bounded mean curvature the multiplicity of $V$ is constant on each
connected component $\Sigma_{\ast}$ of $\Sigma=\reg V\cap(B^{n}_{R}\times(0,\infty))$.

\begin{proof}{\bf{}Proof of Lemma~\ref{mult2}:} If $V$ has multiplicity $N\ge 3$ at a point $(\xi,\tau)\in\Sigma$, then
for suitable $\rho\in (0,\delta)$, $\delta<\tau/2$, we have the conclusions of Lemma~\ref{sheeting} with $L=N\ge 3$. 
Since $|u^{i}_{j}|_{C^{1}}\to 0$ for each $i=1,\ldots,N$ and $\cup_{i=1}^{N}\graph u_{j}^{i}\subset G_{j}$ for each $j$,
for $\sigma\in (0,\rho/2)$ small enough to ensure that $\Sigma$ is a graph of a function with $C^{1}$ norm less than
$\ha$ over the affine tangent hyperplane of $\Sigma$ at the point $(\xi,\tau)$, and for $j$ sufficiently large, we then have
that $U_{j}=\{ x+t\nu(x): x\in \Sigma\cap B_{\sigma}(\xi,\tau)\text{ and } u_{j}^{1}<t<u_{j}^{2}\}$, $\widetilde
U_{j}=\{x+t\nu(x): x\in \Sigma\cap B_{\sigma}(\xi,\tau)\text{ and } u_{j}^{2}<t<u_{j}^{3}\}$ are piecewise $C^{1}$ open
sets.  Also, since $N\ge 3$ and $u_{j}^{1}<u_{j}^{2}<u_{j}^{3}<\cdots <u_{j}^{N}$, we have $G_{j}\cap
(U_{j}\cup\widetilde U_{j})=\emptyset$, so $U_{j},\widetilde U_{j}$ are on opposite sides of $G_{j}$. Since it is then just
a matter of notation, we can, and we shall, assume that $U_{j}$ is contained in the region below the graph $G_{j}$ for
each $j$. For each $j$ let $\nu^{j}$ denote the upward-pointing unit normal function of $G_{j}$, thought of as a
function defined in the whole cylinder $B_{R}^{n}\times(0,\infty)$ which is independent of the $x_{n+1}$ variable.  Thus
$\nu^{j}(x,x_{n+1})=(1+|Du_{j}(x)|^{2})^{-1/2}(-Du_{j}(x),1)$, and the equation~\ref{sme} can be written
$$%
-\op{div}_{\R^{n+1}}\nu^{j}=-\sum_{i=1}^{n+1}D_{i}\nu_{i}^{j}(x,x_{n+1})=H_{j}(x,x_{n+1}),
$$%
valid in the cylinder $B_{R}^{n}\times(0,\infty)$, where $H_{j}$ denotes the mean curvature function defined on
$B_{R}^{n}\times(0,\infty)$ by
$$%
H_{j}(x,x_{n+1})=(m-1)(1+|Du_{j}(x)|^{2})^{-1/2}/u_{j}(x),
$$%
so that $H_{j}(x,x_{n+1})$ is also independent of the $x_{n+1}$ variable.  Observe also that, since $U_{j}$ is contained
in the region below the graph, the outward pointing unit normal of $\partial U_{j}$ on $\graph u_{j}^{1},\graph
u_{j}^{2}$ agrees with $\nu^{j}$ and the remaining part of $\partial U_{j}$ has measure $\le
C\sup(u_{j}^{2}-u_{j}^{1})\to 0$. So applying the divergence theorem over $U_{j}$ gives
$$%
\tint_{U_{j}}H_{j}=\tint_{U_{j}}\op{div}_{\R^{n+1}}\nu^{j} ={\cal{}H}^{n}(\cup_{i=1}^{2}\graph (u_{j}^{i}|\Sigma\cap
B_{\sigma}(\xi,\tau))) + E_{j},
$$%
where $E_{j}\to 0$. On the other hand $H_{j}$ is bounded independent of $j$ on $U_{j}$ (because $U_{j}$ is contained
in the region below the graph $G_{j}$ so at each point $(x,x_{n+1})\in U_{j}$ we have $u_{j}(x)\ge x_{n+1}$ and hence,
since $x_{n+1}>\tau-\rho>\delta$ for any $(x,x_{n+1})\in U_{j}$ and all sufficiently large $j$, we have $u_{j}(x)>\delta$
whenever $(x,x_{n+1})\in U_{j}$, so in particular $H_{j}(x,x_{n+1})\le (m-1)/\delta$ for all sufficiently large $j$ if
$(x,x_{n+1})\in U_{j}$).  Thus the above identity gives ${\cal{}H}^{n}(\graph (u_{j}^{1}|\Sigma\cap
B_{\sigma}(\xi,\tau))\cup \graph (u_{j}^{2}|\Sigma\cap B_{\sigma}(\xi,\tau)))\to 0$, which of course is impossible because
$u_{j}^{i}|\Sigma\cap B_{\sigma}(\xi,\tau)\to 0$ in the $C^{1}$ sense and hence $\liminf_{j\to\infty}{\cal{}H}^{n}(\graph
(u_{j}^{1}|\Sigma\cap B_{\sigma}(\xi,\tau))\cup \graph (u_{j}^{2}|\Sigma\cap B_{\sigma}(\xi,\tau)))\ge
2{\cal{}H}^{n}(\Sigma\cap B_{\sigma}(\xi,\tau))>0$. \end{proof}

\medskip

\begin{state}{\tl{vert}\,Lemma.}%
  Let $\Sigma_{\ast}$ be a connected component of $\Sigma=\reg V\cap (B_{R}^{n}\times (0,\infty))$, equipped with the
same multiplicity as $V$ at each of its points. 

  \emph{(i)} If $\Sigma_{\ast}$ is a vertical cylinder $\Sigma_{0}\times(0,\infty)$, then it has zero mean curvature (i.e.\
minimal) and is stable (i.e.\ it satisfies the stability inequality $\tint_{\Sigma_{\ast}}|A_{\Sigma_{\ast}}|^{2}\zeta^{2}\le
\tint_{\Sigma_{\ast}}|\nabla_{\Sigma_{\ast}}\zeta|^{2}$ for each $\zeta\in C^{1}_{c}(B_{R}^{n}\times(0,\infty))$), and
furthermore $\overline{\Sigma}_{\ast}\cap \overline{\Sigma\setminus\Sigma}_{\ast}=\emptyset$ (where $\overline S$
means closure of $S$ in $B_{R}^{n}\times\R$), and $\Sigma_{\ast}$ has multiplicity~$2$.

\emph{(ii)} If $\Sigma_{\ast}$ has multiplicity~$2$ then it is a vertical cylinder $\Sigma_{0}\times(0,\infty)$ 
  \end{state}

  \begin{proof}{\bf{}Proof of~(i):} The mean curvature of $\Sigma$ is $\le C/\sigma$ in the region
$x^{n+1}>\sigma$; but of course since $\Sigma$ is a vertical cylinder this shows that $\Sigma$ is minimal (i.e.\ mean
curvature zero).

To prove that $\Sigma$ is stable in $\R^{n}\times(0,\infty)$ (or equivalently $\Sigma_{0}$ is stable in $\R^{n}$) we let
$\zeta\in C^{1}_{c}(\R^{n}\times(0,\infty))$ be arbitrary and for any given $K$ let $\zeta_{K}(x,x^{n+1})=
\zeta(x,x^{n+1}-K)$, so that support $\zeta_{K}\subset (K,\infty)$. Then by the inequality \ref{ssf1} we have
$$%
\tint_{\Sigma} |\nabla_{\Sigma}\zeta_{K}|^{2}  \le (1-\epsilon)^{-1}\tint_{\Sigma}|A_{\Sigma}|^{2}\zeta_{K}^{2}
+\epsilon^{-2}K^{-2}\tint_{\Sigma}\zeta_{K}^{2} 
$$%
Since $\Sigma$ is cylindrical this can be written 
$$%
\tint_{\Sigma} |\nabla_{\Sigma}\zeta|^{2}  \le (1-\epsilon)^{-1}\tint_{\Sigma}|A_{\Sigma}|^{2}\zeta^{2}
+\epsilon^{-2}K^{-2}\tint_{\Sigma}\zeta^{2} 
$$%
so by first letting $K\to \infty$ and then letting $\epsilon\downarrow0$ we obtain the stability inequality
as claimed.

Let $S(\Sigma)=\{(x,\xi)\in B_{R}^{n}\times(\R^{m}\setminus\{0\}):(x,|\xi|)\in\Sigma\}$, interpreted as an integer
multiplicity varifold equipped at each point $(x,\xi)\in S(\Sigma)$ with the multiplicity of $\Sigma$ at the point
$(x,|\xi|)$. Then $S(\Sigma)$ is the varifold limit of $S(u_{j})$ for some sequence $u_{j}$ of singular solutions $u_{j}$
of~\ref{sme} on $B_{R}^{n}$ and hence $S(\Sigma)$ is a stationary integer multiplicity varifold in
$B_{R}^{n}\times\R^{m}$. Notice that, since $m\ge 2$, $\R^{n}\times 0$ has ${\cal{}H}^{n+m-1}$-dimensional measure
zero, so in terms of varifold convergence $S(\Sigma)=\lim S(u_{j})$, and each $S(u_{j})$ is stationary in
$B_{R}^{n}\times\R^{m}$, and hence $S(\Sigma)$ is indeed stationary in $B_{R}^{n}\times\R^{m}$ rather than merely
(locally) stationary in $B_{R}^{n}\times(\R^{m}\setminus\{0\})$, although a-priori we need to allow the possibility that
$\sing S(\Sigma)(=\overline{S(\Sigma)}\setminus S(\Sigma))$ might include a subset of $B_{R}^{n}\times\{0\}$, even
perhaps a subset of positive $n$-dimensional Hausdorff measure; ultimately we shall prove in~\S\ref{dimsing} that at
most a set of Hausdorff dimension $n-2$ occurs here.

Next observe that we can apply the maximum principle of Ilmanen \cite{Il96} to assert $\overline{S(\Sigma_{\ast})}\cap
\overline{S(\Sigma)\setminus S(\Sigma_{\ast})}=\emptyset$ as claimed; notice that, since we are working with $N=n+m-1$
dimensional hypersurfaces, to literally apply \cite{Il96} we need to check that $\overline{S(\Sigma)\setminus
S(\Sigma_{\ast})}\cap \overline{S(\Sigma_{\ast})}$ has $N-2\,(=n+m-3)$-dimensional Hausdorff measure zero, but that
hypothesis was only used in \cite{Il96} to justify application of the regularity theory \cite{SS81}, and that really only uses
that the ``\emph{capacity}'' of any compact $K\subset \overline{S(\Sigma)}$ is zero,  meaning that there is a sequence
$\zeta_{j}$ of $C_{c}^{1}$ functions with $\zeta_{j}=0$ in a neighborhood of $K$, $\zeta_{j}\equiv 1$ in $\{x:\dist(x,K)\ge
1/j\}$ and $\tint_{S(\Sigma)} |\nabla \zeta_{j}|^{2}\to 0$.  Thus it suffices to check that $\overline{S(\Sigma)\setminus
S(\Sigma_{\ast})}\cap \overline{S(\Sigma_{\ast})}$ has \emph{locally finite} (rather than zero) $(N-2)$-dimensional
Hausdorff measure, because a compact set $K$ with finite $(N-2)$-dimensional Hausdorff measure automatically has
capacity zero. (See e.g.\ \cite[Theorem 4.16]{EG}.)  In the present case we have $\overline{S(\Sigma_{\ast})}\cap
(B_{R}^{n}\times\{0\})$ is a set of locally finite $n-1$ dimensional Hausdorff measure, and $n-1=(n+m-1)-m=N-m\le
N-2$, so indeed \cite{Il96} can be applied, and $\overline{S(\Sigma)\setminus S(\Sigma_{\ast})}\cap
\overline{S(\Sigma_{\ast})}=\emptyset$ as claimed. We should mention here that in any case recent work of
Wickramasekera (\cite{Wic14}) shows that the application of the maximum principle requires only that the
$(N-1)$-dimensional Hausdorff measure of the intersection of the supports is zero and of course we have that here. So
technically we did not need to include the discussion concerning capacity in relation to Hausdorff $(N-2)$ dimensional
measure, but it was included for the reader's convenience since the work of Wickramasekera is lengthy and deep.

Finally we have to check that $\Sigma_{\ast}=\Sigma_{0}\times(0,\infty)$ has multiplicity~$2$.  By Lemma~\ref{mult2} the
multiplicity is either $1$ or $2$.  Of course since the multiplicity is constant on the regular set it suffices to assume $0\in
\Sigma_{0}\times\{0\}$ and show that it is not possible for $\Sigma_{\ast}\cap (B_{\sigma}^{n}\times (0,1))$ to have
multiplicity~$1$.  In view of what has already been established above, we can choose $\sigma>0$ small enough to ensure
that $B_{\sigma}^{n}\times(0,1)\cap (S(\Sigma\setminus\Sigma_{\ast}))=\emptyset$.  Let $u_{j}$ be as above, so we have
varifold convergence $\Sigma_{\ast}=\lim G_{j}$ in $B_{\sigma}^{n}\times(-1,1)$, where $G_{j}=G(u_{j})$ considered as a
multiplicity~1 varifold.  We can of course use the orientation of $G_{j}$ given by the upward pointing unit normal and
view $G_{j}$ as a multiplicity~1 current in $B_{\sigma}\times(-1,1)$ with $\partial G_{j}=0$. The varifold convergence
$G_{j}\to\Sigma_{\ast}$ is multiplicity~1, so using the local $C^{k}$ convergence guaranteed by Lemma~\ref{sheeting}, we
can appropriately orient $\Sigma_{\ast}$ so that the convergence of $G_{j}$ to $\Sigma_{\ast}$ in
$B_{\sigma}\times(-1,1)$ is also in the weak sense of currents (thus $\tint_{G_{j}}\omega\to \tint_{\Sigma_{\ast}}\omega$
for each fixed smooth $n$-form $\omega$ with compact support in $B_{\sigma}^{n}\times(-1,1)$), and hence
$\partial\Sigma_{\ast}=\lim \partial G_{j} =0$.  But of course $\partial\Sigma_{\ast}=\pm \Sigma_{0}\times \{0\}$, a
contradiction.
\end{proof}

\begin{proof}{\bf{}Proof of~(ii):} Take $(\xi,\tau)\in\Sigma_{\ast}$ and let $u_{j}^{1}<u_{j}^{2}$ be as in the sheeting
lemma~\ref{sheeting}, so that $|u^{i}_{j}|_{C^{2}}\to 0$ for $j=1,2$. Now if $\Sigma_{\ast}$ has normal $\nu$ with
$\nu_{n+1}(\xi,\tau)>0$ (i.e.\ $\Sigma_{\ast}$ has a non-vertical tangent hyperplane at $(\xi,\tau)$) then the $C^{1}$
convergence of $u_{j}^{i},u_{j}^{2}$ ensures that the vertical line $\{\xi\}\times(0,\infty)$ intersects both $\graph
u_{j}^{1}$ and $\graph u_{j}^{2}$ for sufficiently large $j$, which contradicts the fact that $G_{j}$ is a graph over
$B_{R}^{n}$.  Thus $\nu_{n+1}\equiv 0$ and in particular the mean curvature of $\Sigma_{\ast}$ (which is
$(m-1)\nu_{n+1}/x_{n+1}$) is identically zero.  Thus $\Sigma_{\ast}$ is vertical (i.e.\ $\nu_{n+1}\equiv 0$ on
$\Sigma_{\ast}$) and has zero mean curvature.

To show that $\Sigma_{\ast}$ is actually a cylinder $\Sigma_{0}\times(0,\infty)$ we have simply to prove that
$(\xi,\tau)\in\Sigma_{\ast}\Rightarrow$ the whole ray $\{\xi\}\times(0,\infty)\subset\Sigma_{\ast}$. Since
${\cal{}H}^{n-1}(\sing V\cap (B^{n}_{R}\times(0,\infty)))=0$ we must have
${\cal{}H}^{n-1}(P(\overline{\Sigma}_{\ast}\setminus\Sigma_{\ast})\cap (B_{R}\times\{0\}))=0$, where $P$ is the orthogonal
projection of $\R^{n}\times\R^{m}$ onto $\R^{n}\times\{0\}$.  With $S=\overline{\Sigma}_{\ast}\setminus\Sigma_{\ast}$
we then have ${\cal{}H}^{n}\bigl(P(S)\times (0,\infty)\bigr)=0$, and since $\Sigma_{\ast}$ is smooth embedded with
$\nu_{n+1}\equiv 0$ it is clear that $\Sigma_{\ast}\setminus \bigl(P(S)\times (0,\infty)\bigr)$ is a union of vertical rays
$\{\xi\}\times(0,\infty)$ and hence is a cylinder which agrees up to a set of ${\cal{}H}^{n}$-measure zero with
$\Sigma_{\ast}$, and in particular is dense in $\Sigma_{\ast}$.  Thus $\overline{\Sigma}_{\ast}$ is the closure (in
$B_{R}^{n}\times(0,\infty)$) of a cylinder and hence is a cylinder.
\end{proof}

\section{Tangent Cones of Singular Solutions}\label{tcs}

Let $u$ be an arbitrary singular solution of \ref{sme} in $B_{R}^{n}$ and suppose $0\in \sing u$. Since $S(u)$
is a stationary integer multiplicity varifold we can take tangent cones of $S(u)$ at $0$. Thus for each sequence
$\lambda_{j}\downarrow 0$ we can take a subsequence (still denoted $\lambda_{j}$) such that
$\lambda_{j}^{-1}S(u)$ converges in the varifold sense in $\R^{n}\times \R^{m}$ to a cone $\mathbb{S}$ (i.e.\
$h\mathbb{S}=\mathbb{S}$ for each $h>0$), and by construction $\mathbb{S}$ is invariant under all orthogonal
transformations of $\R^{n}\times\R^{m}$ which leave the first $n$ coordinates fixed. 

In terms of $u$ this means that, with $u_{j}(x)=\lambda_{j}^{-1}u(\lambda_{j}x)$ for $x\in
\smash{B_{R/\lambda_{j}}^{n}}$ and $G_{j}=G(u_{j})$ (the graph of $u_{j}$), $G_{j}$ converges in the varifold sense to
an integer multiplicity cone $\C$ in $\R^{n}\times\R$.  Such $\C$ is called a tangent cone to $G(u)$ at $0$ and of
course $\C$ is indeed a cone in the sense that $h\C=\C$ for each $h>0$.
  
\begin{state}{\tl{tan-cones} Lemma.} If $u$ is a singular solution of \emph{\ref{sme}} on the ball $B_{R}^{n}$, and if $0\in
\sing u$, then any tangent cone $\C$ of $G(u)$ at $0$ (obtained as described above) is of multiplicity~$1$.
  \end{state}

\begin{proof}{\bf Proof:}  Let $u$ be any singular solution of~\ref{sme} on $B_{R}^{n}$ with
$0\in\sing u$, and, with the terminology introduced above,  let $\C$ be any tangent cone to $G(u)$ at~$0$.

By the regularity discussion of \S\ref{sec-regularity}, in the region $x^{n+1}>0$ $\sing\C$ has Hausdorff dimension $\le
n-7$ for $n\ge 7$, empty for $n\le 6$, discrete for $n=7$.  Since $h\reg\C=\reg\C$ for each $h>0$, each
connected component of $\reg\C$ has zero in its closure taken in $\R^{n}\times\R$.  Thus if $\C\cap
(\R^{n}\times(0,\infty))$ has a vertical component $\Sigma_{0}\times(0,\infty)$ then by Lemma~\ref{vert} it has
multiplicity~2 and is all of $\C\cap (\R^{n}\times(0,\infty))$. We claim that in this case \emph{every} tangent cone would
also have to be vertical. Indeed the set of all tangent cones of $G(u)$ at $0$ is evidently connected, so if there is a
non-vertical tangent cone of $G(u)$ at $0$ then there would have to be a sequence of $\C_{j}$ of non-vertical tangent
cones converging to a vertical tangent cone $\C$ in the varifold sense. Also, again using Lemma~\ref{vert}, each $\C_{j}$
must have multiplicity~1, so using standard Harnack theory locally in $\Sigma=\reg \C\cap (\R^{n}\times(0,\infty))$ and
the sheeting lemma~\ref{sheeting} in the manner of \cite{Il96} we would conclude that there is a positive smooth solution
$v$ of the Jacobi field equation $\Delta_{\Sigma}v+|A_{\Sigma}|^{2}v=0$ on $\Sigma$. Since $\Sigma$ is obtained by a
rescaling of the difference in heights of the ``sheets'' (i.e.\ the difference $u^{j}_{1}-u^{j}_{2}$ in the terminology of
Lemma~\ref{sheeting}), and each $\C_{j}$ is a cone, we know that $v$ is homogeneous of degree 1 on $\Sigma$. Thus
$v$ is homogeneous degree 1 superharmonic (i.e.\ $\Delta_{\Sigma} v\le 0$) on $\Sigma$. Then of course $\min\{v,K\}$
is a bounded weakly superharmonic on $\Sigma$ for each constant $K>0$. According the argument of \cite{Il96}, since
$\C$ is a cone we can then use the mean value inequality for $\min\{v,K\}$:
$$%
\tint_{\Sigma\cap B_{1}^{n+1}} \min\{v,K\} \le  \rho^{-n}\tint_{\Sigma\cap B_{\rho}^{n+1}} \min\{v,K\} ,
\quad \forall\,\rho\in(0,1),\, K>0. 
$$%
Since $\Sigma$ is a cone and since $v$ is homogeneous of degree~1 we can change variable in the integral on the right
side to give
$$%
\rho^{-n}\tint_{\Sigma\cap B_{\rho}^{n+1}} \min\{v,K\}=\tint_{\Sigma\cap B_{1}^{n+1}}\min\{\rho v,K\} %
$$%
and  so the above inequality gives 
$$%
\tint_{\Sigma\cap B_{1}^{n+1}} (\min\{v,K\}-\min\{\rho v,K\}) \le 0. 
$$%
But of course $\min\{\rho v,K\}\le \min\{v,K\}$ so the integrand is non-negative here and hence
$$%
\min\{v,K\}=\min\{\rho v,K\}  \text{ ${\cal{}H}^{n}$-a.e.\ on $\Sigma$},
$$%
which of course is impossible for sufficiently large $K$. 

So we conclude that either \emph{all} the tangent cones of $G(u)$ at $0$ are vertical (hence by Lemma~\ref{sheeting}
also have connected multiplicity~2 regular set), or else there are no tangent cones which are vertical. We claim that the
former case, when all tangent cones of $G(u)$ at $0$ are vertical, multiplicity~2, and with connected regular set cannot
occur.  Indeed in this case we have that each tangent cone $\mathbb{S}$ of $S(u)$ at $0$ corresponds to such a
cylindrical tangent cone $\C$ of $G(u)$ at $0$ according to $\mathbb{S}=\{(x,\xi):(x,|\xi|)\in \C\}$. Thus, with
$N=n+m-1$, in particular every tangent cone $\mathbb{S}$ of $S(u)$ at $0$ has connected regular set and a singular set
(in $\R^{N+1}$) of locally finite ${\cal{}H}^{N-2}$-measure. But, according to   \cite[Theorem B]{Il96} in these
circumstances tangent cones of multiplicity~1 minimal submanifolds  must in fact also be multiplicity~1.  To be strictly
correct we should mention that Ilmanen actually proved his Theorem~B only for stable multiplicity~1 minimal
hypersurfaces with singular sets of ${\cal{}H}^{N-2}$ measure zero, but, as explained in the discussion in the proof of
Lemma~\ref{vert}, Ilmanen's argument applies in the present setting where we have the sheeting lemma~\ref{sheeting}
and the bound ${\cal{}H}^{N-2}(K)<\infty$ for each tangent cone $\mathbb{S}$ of $S(u)$ at $0$ and each compact subset
$K\subset \sing\mathbb{S}$.  \end{proof}

\section[Multiplicity 1 Cones in $\overline{\cal{}G}_{\infty}$]{Multiplicity 1 Cones in
\boldmath{$\overline{\cal{}G}_{\infty}$}}\label{harnack-sec}

\begin{state}{\tl{closure} Lemma.}  Let ${\cal{}C}$ be the family of all multiplicity~$1$ cones $\C\in
\overline{\cal{}G}_{\infty}$ (thus each $\C\in{\cal{}C}$ is a multiplicity~$1$ cone---$\C$ has multiplicity~$1$ at each
point of $\reg\C$ and $\lambda\C=\C$ for each $\lambda>0$---and there exists a sequence $u_{j}$ of solutions
of~\emph{\ref{sme}} on $B_{R_{j}}^{n}$ with $R_{j}\to\infty$ and $G(u_{j})$ converging $\C$ in the varifold sense in
$\R^{n+1}$).

Then ${\cal{}C}$ is a compact subset of $\overline{\cal{}G}_{\infty}$; thus for each sequence $\C_{j}\subset {\cal{}C}$
there is $\C\in{\cal{}C}$ and a subsequence $\C_{j'}$ converging to $\C$ in the varifold sense.
\end{state}

\begin{proof}{\bf{}Proof:} In view of the volume bounds \ref{vol-bds} the Allard compactness theorem guarantees that for
any sequence $\C_{j}\subset {\cal{}C}$ there is a subsequence $\C_{j'}$ converging to $\C\in\overline{\cal{}G}_{\infty}$ in
the varifold sense.

So to complete the proof of compactness we just have to prove that ${\cal{}C}$ is closed in $\overline{\cal{}G}_{\infty}$;
i.e.\ if $\C_{j}\in{\cal{}C}$ with $\C_{j}\to \C$ in the varifold sense, then $\C$ is a multiplicity~$1$ cone in
$\overline{\cal{}G}_{\infty}$.  Certainly $\C_{j}\to \C$ implies that $\C$ is in $\overline{\cal{}G}_{\infty}$ (because
$\overline{\cal{}G}_{\infty}$ is closed by definition), and trivially $\C$ is a cone, so we just have to check that $\C$ has
multiplicity~$1$.

Otherwise, using the sheeting lemma \ref{sheeting} exactly as in the proof of Lemma~\ref{tan-cones} we would have a
component $\Sigma$ of $\reg\C$ and homogeneous degree one smooth positive $v$ on $\Sigma$ with
$\Delta_{\Sigma}v\le 0$, and by the same argument as in the proof of Lemma~\ref{tan-cones}, using the mean value
inequality on $\Sigma$, this is impossible.
\end{proof}

The following lemma establishes a Harnack theory for certain supersolutions  on domains in $\C$, $\C\in {\cal{}C}$.

\begin{state}{\tl{harnack} Lemma (Harnack for Supersolutions.)}
There is $\lambda=\lambda(m,n)\in (0,\frac{1}{8}]$ such that if  $\beta>0$,  $\C\in {\cal{}C}$  and
if $v$ is a bounded positive $C^{1}$ function on $\Sigma\cap B^{n+1}_{1/2}(e_{n+1})$ ($\Sigma=\reg\C$) which
satisfies an inequality of the form
$$%
\Delta_{\Sigma}v+b\cdot\nabla_{\Sigma}v +c v\le 0 \text{ on }\Sigma\cap B^{n+1}_{1/2}(e_{n+1}),  
$$%
with $|b|,|c|\le \beta$, then 
$$%
\tint_{\Sigma\cap B_{\lambda}^{n+1}(e_{n+1})}v \le C \inf_{\Sigma\cap B^{n+1}_{\lambda}(e_{n+1})} v, \quad
C=C(m,n,\beta).
$$%
\end{state}

{\bf{}Remark:} Notice that this is uniformly applicable over all $\C\in{\cal{}C}$, because the constant $C$ depends only
on $m,n,\beta$ and not on the particular $\C\in{\cal{}C}$.

\begin{proof}{\bf{}Proof of Lemma~\ref{harnack}:} Since we have a suitable Sobolev inequality (see \cite{MicS73}) it is
well known (see e.g.\ the discussion in \cite{BomG72}) that one can apply a Harnack theory in $\Sigma\cap
B_{1}^{n+1}(e_{n+1})$, $\Sigma=\reg\C$, for positive bounded $v$ satisfying
$\Delta_{\Sigma}v+b\cdot\nabla_{\Sigma}v+cv\le 0$ in $B_{1}^{n+1}(e_{n+1})$, provided $|b|,|c|\le \beta$ and provided we
have a suitable Poincar\'e inequality
$$%
\min_{\mu\in \R}\tint_{\Sigma\cap B_{\lambda}^{n+1}(e_{n+1})}|h-\mu|\le
\gamma\tint_{\Sigma\cap B^{n+1}_{1/4}(e_{n+1})}|\nabla_{\Sigma} h|
\pdl{an-poinc}
$$%
for $h\in C^{1}(\Sigma)$.  As is well known, such an inequality is implied by the geometric inequality
$$%
\min\{{\cal{}H}^{n}(E\cap B^{n+1}_{\lambda}(e_{n+1})),{\cal{}H}^{n}(\Sigma\cap B^{n+1}_{\lambda}(e_{n+1})\setminus
E)\} \le \gamma{\cal{}H}^{n-1}(\Gamma\cap B^{n+1}_{1/4}(e_{n+1})), %
\pdl{poinc}
$$%
for sets $E\subset \Sigma\cap B_{1/4}^{n+1}(e_{n+1})$ where $\Gamma=(\overline{\!E}\setminus E)\cap
B^{n+1}_{1/4}(e_{n+1})$. (The inequality \ref{an-poinc} follows from \ref{poinc} by taking $E=E_{t}=\{x:h(x)<t\}$ and then
integrating with respect to $t$ and applying the coarea formula.)

Using the Sobolev inequality \cite{MicS73} and the volume bounds of \ref{vol-bds}, it is a standard fact that there
is $\eta=\eta(m,n)>0$ such that the inequality~\ref{poinc} holds with $\gamma=\gamma(m,n)$ if
$$%
\min\{{\cal{}H}^{n}(E\cap B^{n+1}_{2\lambda}(e_{n+1})),  %
       {\cal{}H}^{n}(\Sigma\cap B^{n+1}_{2\lambda}(e_{n+1})\setminus E)\} \le \eta\lambda^{n}. %
\pdl{small}
$$%
So we only have to prove~\ref{poinc} subject to the extra assumption that 
$$%
\min\{{\cal{}H}^{n}(E\cap B^{n+1}_{2\lambda}(e_{n+1})),{\cal{}H}^{n}(\Sigma\cap B^{n+1}_{2\lambda}(e_{n+1})\setminus
E)\} \ge \eta\lambda^{n}. %
\pdl{large}
$$%
As a preliminary to this we first claim that there is $\lambda_{0}=\lambda_{0}(m,n)\in (0,\frac{1}{8}]$ such that 
if $\C\in{\cal{}C}$ and if $\Sigma$ is a connected component of $\reg \C$ then
\ab{3pt}{3pt}{%
\begin{align*}%
  &\ptg{lambda}\text{$\Sigma^{1},\Sigma^{2}$   %
distinct connected components of $\Sigma\cap B^{n+1}_{1/4}(e_{n+1})\Longrightarrow$} \\
\noalign{\vskip-3pt}
  &\hskip1.2in  \text{either  $\Sigma^{1}\cap B^{n+1}_{\lambda_{0}}(e_{n+1})=\emptyset$ or $\Sigma^{2}\cap
B^{n+1}_{\lambda_{0}}(e_{n+1})=\emptyset$. }
\end{align*}}%
Indeed otherwise there would be a sequence $\C_{j}\in {\cal{}C}$ and components $\Sigma_{j}^{1},\Sigma_{j}^{2}$ of
$\Sigma_{j}\cap B^{n+1}_{1/4}(e_{n+1})$ ($\Sigma_{j}=\reg\C_{j}$) such that $\Sigma_{j}^{i}\cap
\smash{B^{n+1}_{1/j}(e_{n+1})}\neq\emptyset$ for $i=1$ and $i=2$. But then by the volume bounds~\ref{vol-bds} and
the Allard compactness theorem, $\Sigma_{j}^{i}$ converges in the varifold sense to a bounded mean curvature
integer multiplicity varifold $V^{i}$ in $B^{n+1}_{1/4}(e_{n+1})$ and $e_{n+1}\in\spt V^{1}\cap \spt V^{2}$ and of
course $\spt V^{i}\subset\spt C$ for $i=1,2$.  Then by the maximum principle \cite{Il96} we must have
${\cal{}H}^{n-2}(\spt V^{1}\cap \spt V^{2})\neq 0$.  But,  since ${\cal{}H}^{n-2}(\sing \C\cap \{x^{n+1}>0\})=0$, we
must then have at least one  point $\xi\in\reg\C\cap\spt V^{1}\cap\spt V^{2}$. Since $\reg\C$ is smooth the constancy
theorem implies that both $V^{1}$ and $V^{2}$ are positive integer multiples of $\C$ in a neighborhood of $\xi$
and by construction $V^{1}+V^{2}\le \C$ in $B^{n+1}_{1/4}(e_{n+1})$, so $\C$ would have multiplicity $\ge 2$
in a neighborhood of $\xi$, which contradicts Lemma~\ref{closure}.

Thus there is indeed a $\lambda_{0}=\lambda_{0}(m,n)$ as in~\ref{lambda}.  With this $\lambda_{0}$ we claim that
\ref{poinc} holds with $\lambda=\lambda_{0}/3$ and for
some $\gamma=\gamma(m,n)$. 

If \ref{poinc} is false, then there would be a sequence $\C_{j}\subset {\cal{}C}$ and measurable subsets $E_{j}\subset
\Sigma_{j}\cap \smash{B_{1/4}^{n+1}(e_{n+1})}$, $\Sigma_{j}=\reg\C_{j}$, such that \ref{large} holds with $E=E_{j}$ and
${\cal{}H}^{n}(\Gamma_{j})\to 0$, where $\Gamma_{j}=(\overline{\!E}_{j}\setminus E_{j})\cap
\smash{B^{n+1}_{1/4}(e_{n+1})}$.  Let $F_{j}=\Sigma\cap \smash{B^{n+1}_{1/4}(e_{n+1})}\setminus E_{j}$ and view
$E_{j},F_{j}$ as varifolds in $\smash{B^{n+1}_{1/4}(e_{n+1})}$.  Since the mean curvature of $C_{j}\le 2(m-1)$ in
$B_{1/4}(e_{n+1})$ the first variations $|\delta E_{j}(\zeta)|,|\delta F_{j}(\zeta)|$ are evidently $\le
2(m-1)\|\zeta\|_{L^{1}}+{\cal{}H}^{n-1}(\partial\Gamma_{j})\sup|\zeta|\to 2(m-1)\|\zeta\|_{L^{1}}$ for $\zeta\in
C^{1}_{c}(B_{1/4}(e_{n+1});\R^{n+m})$ (viewing $E_{j},F_{j}$ as multiplicity~1 varifolds), so by the Allard compactness
theorem $E_{j},F_{j}$ converge to bounded mean curvature integer multiplicity varifolds $V^{1},V^{2}$ in
$B^{n+1}_{1/4}(e_{n+1})$ with both $\spt V^{1}\cap B^{n+1}_{2\lambda}(e_{n+1})\neq \emptyset$ and $\spt V^{2}\cap
B^{n+1}_{2\lambda}(e_{n+1})\neq \emptyset$ and $\spt V^{i}\cap B^{n+1}_{1/4}(e_{n+1})\subset \spt \C \cap
B^{n+1}_{1/4}(e_{n+1})$ for $i=1,2$.  By the constancy theorem $V^{i}$ is a sum of positive integer multiplicities
of some components of $\reg \C\cap B^{n+1}_{1/4}(e_{n+1})$ for $i=1,2$ and for each $i=1,2$ one of these
components must intersect $B^{n+1}_{2\lambda}(e_{n+1})$.

But since $2\lambda=2\lambda_{0}/3<\lambda_{0}$, we then conclude that these two components must coincide so again
$\C$ has a connected component of multiplicity~2 which again contradicts Lemma~\ref{closure}.

Thus we do have~\ref{poinc}, hence \ref{an-poinc}, and so by the discussion of \cite{BomG72} we do indeed have the
relevant Harnack theory for supersolutions, so the inequality for $v$ claimed in the lemma is proved.~\end{proof}

\begin{state}{\tl{bd-slope} Theorem.}
  If ${\cal{}C}$ is as above then there is a bound, depending on $m,n$ only, on the slope of rays for the cones
$\C\in{\cal{}C}$; in fact there is a $K_{0}=K_{0}(m,n)>0$ such that if $\C\in{\cal{}C}$ then there is a
$C^{0,{1/2}}(\R^{n})$ homogeneous degree $1$ singular solution  $\varphi$ of~\emph{\ref{sme}} with
$G(\varphi)=\C$ and $\varphi(x)\le K_{0}|x|$ for all $x\in\R^{n}$.
\end{state}

\begin{proof}{\bf{}Proof:} If the first claim is false then there is a sequence $\C_{j}\in {\cal{}C}$ such that, for each
$\sigma>0$, $\C_{j}\cap (B^{n}_{\sigma}(0)\times (1,\infty))\neq \emptyset$ for all $j$ sufficiently large (depending on
$\sigma$), and by Lemma~\ref{closure} there is a subsequence of $\C_{j}$ (still denoted $\C_{j}$) and $\C\in{\cal{}C}$
with $\C_{j}\to \C$. By construction $\spt\C\cap (B^{n}_{\sigma}\times(1,\infty))\neq \emptyset$ for each $\sigma>0$, so
$\{(0,t):t>0\}\subset\spt\C$.  Also on $\reg \C$ we have by~\ref{limits2} that
$$%
\Delta_{\C}\nu_{n+1} + |A_{\C}|^{2}\nu_{n+1} = -(m-1)e_{n+1}\cdot \nabla_{\C}((x^{n+1})^{-1}\nu_{n+1}),
$$%
and hence 
$$%
\Delta_{\C}\nu_{n+1} +(m-1)e_{n+1}\cdot \nabla_{\C}((x^{n+1})^{-1}\nu_{n+1}) \le 0.
$$%
But we can then apply Lemma~\ref{harnack}, hence, with $\Sigma$ any connected component of $\reg\C\cap
(\R^{n}\times(0,\infty))$,
$$%
\tint_{\Sigma\cap B^{n+1}_{\lambda}(e_{n+1})}\nu_{n+1} \le C\inf_{\Sigma\cap B^{n+1}_{\lambda}(e_{n+1})}\nu_{n+1}.
$$%
In view of the inclusion $\{(0,t):t>0\}\subset\spt\C$ the right side here is zero, so $\nu_{n+1}$ is identically zero on
$\Sigma\cap B_{\lambda}^{n+1}(e_{n+1})$, hence on all of $\Sigma$.  Thus $\reg\C \cap (\R^{n}\times(0,\infty))$
contains a vertical cylinder $\Sigma_{0}\times(0,\infty)$. But by Lemma~\ref{vert} such a vertical cylinder has
multiplicity~2, contradicting Lemma~\ref{closure}.

Now let $\C$ be any tangent cone of $G(u)$. Then there is a sequence $\lambda_{j}\downarrow 0$ and $G(u_{j})\to\C$
on $B_{1}(0)\times(0,\infty)$, where $u_{j}(x)=\lambda_{j}^{-1}u(\lambda_{j}x)$.  In view of the uniform bound (by
$K_{0}$ say) on the  slope of the rays of $\C$ established above, we have $\spt\C\cap (B_{1}^{n}\times\R) \subset
B^{n}_{1}\times [0,K_{0}]$,  so by Lemma~\ref{holdercor} the $u_{j}$ are uniformly bounded in
$C^{0,{\frac{1}{2}}}(B_{\!\frac{1}{2}})$ and hence there is a subsequence (still denoted
$u_{j}$) converging uniformly to  $\varphi$ on $B_{\!\frac{1}{2}}$ (and of course the convergence
is locally in the $\smash{C^{k}}$ sense for each $k$ on the open set where $\varphi>0$), so $\spt \C\cap
(B_{\!\frac{1}{2}}\times(0,\infty))= \reg \C\cap (B_{\!\frac{1}{2}}\times(0,\infty)) 
=\graph \varphi|\{x:\varphi>0\}$, and hence $|D_{r}\varphi|\le K_{0}$ on $\{x:\varphi>0\}$ ($D_{r}\varphi=|x|^{-1}x\cdot
D\varphi$). Evidently then $\varphi(x)\le K_{0}|x|$ for all $x\in B_{\!\frac{1}{2}}$.  But since $\C$ is a
cone $\varphi$ is the restriction to $B_{\!\frac{1}{2}}$ of a homogeneous degree~1 function
on $\R^{n}$.  This completes the proof of Theorem~\ref{bd-slope}.~\end{proof}

\section{Gradient Estimates} \label{grad-ests}

We can now prove the gradient bounds mentioned in \S\hskip-1pt\ref{intro}:

\begin{state}{\tl{grad0}\,Theorem.}%
  Let $u$ be a regular or singular solution of the SME (i.e.\ \emph{\ref{sme}}) on a ball $B_{\rho}(y)$ with
$\sup_{B_{\rho}(y)} u < M$. Then, for any $\theta \in (0,1)$,
$$%
\sup_{B_{\theta\rho}(y)}|D u| < C,
$$%
where $C$ depends only on $m$, $n$, $M/\rho$, and $\theta$.
\end{state}

\begin{proof}{\bf{}Proof:}
  By scaling it suffices to prove the theorem with $\rho=1$.  If this is false, there would exist constants $M>0$, $\theta\in
(0,1)$, a sequence $\{x_j\} \in B_{\theta}(y_{j})$, and a corresponding sequence $\{u_j\}$ of (possibly singular) solutions
of \ref{sme} defined on $B_{1}(y_{j})$ with $\sup u_{j} < M$ such that $u_{j}(x_{j})>0$ and $|D u_j (x_j)| \uparrow
\infty$. Because the solutions of \ref{sme} are invariant under translations of $\R^n$, translating the $\graph u_j$ by
$(x_j,0)$, we can assume that $x_j=0$ for all $j$.  Thus we have $u_{j}$ defined at least over $B_{R}(0)$, where
$R=1-\theta$, and $\sup_{B_{R}(0)} u_{j}\le M$, $|Du_{j}(0)|\to \infty$.  But by virtue of the fixed bounds on the gradient
of $u_{j}^2$ (Corollary \ref{holdercor}), we have
\ab{3pt}{3pt}{%
$$%
\sup_{j}u_{j}(0) |D u_{j}(0)| < \infty,
$$}%
and hence
\ab{3pt}{3pt}{%
$$%
u_{j}(0) \rightarrow 0. %
\pdl{ujto0}
$$}%

Let $K >M/(1-\theta)$ be a fixed constant (to be chosen later), and $W$ be the open circular cone
\ab{0pt}{3pt}{%
\begin{align*}%
  W&= \{(\xi, \tau) \in \R^n \times (0,\infty): \tau> K|\xi|\} \\
\noalign{so that the boundary of $W$ is a union of rays $\tau=K|\xi|$ of slope $K$, and let \vskip3pt}
  W_{t}&= W\cap (\R^{n}\times\{t\}).
\end{align*}}%
Note that $(0,u_j(0)) \in W$, but $\sup_{B_R} u_{j} \le M$ and $K>M/(1-\theta)$, so $G_j \cap W_{t}=\emptyset$ for
$t= M$; thus $G_{j}$ intersects $W$ non-trivially at some heights $t\ge u_{j}(0)$, but eventually leaves $W$ completely. 
We let
\ab{3pt}{3pt}{%
$$%
h_{j} =\inf\{t:t>u_{j}(0) \text{ and }G_{j}\cap W_{t}=\emptyset \},
$$}%
so then
\ab{3pt}{3pt}{%
$$%
G_{j}\cap W_{t}\neq \emptyset \text{ for each }t\in [u_{j}(0),h_{j}) \text{ and }\sup_{B_{h_{j}/K}(0)}u=h_{j}. %
\pdl{nonempty}
$$}%
We claim that
\ab{3pt}{3pt}{%
$$%
\frac{h_j}{u_j(0)} \rightarrow \infty\text{ as }j \rightarrow \infty.
\pdl{hoveru}
$$}%
To check this, rescale to give $\widetilde{G}_{j}$ according to
\ab{3pt}{3pt}{%
$$%
\widetilde{G}_j=(u_{j}(0))^{-1} G_j
$$}%
Then $\widetilde{G}_{j}=\graph \widetilde{u}_{j}$, where $\widetilde{u}_{j}=\lambda_{j}^{-1}u_{j}(\lambda_{j}x)$, with
$\lambda_{j}=u_{j}(0)$, is a (possibly singular) solution of \ref{sme}, $|D\widetilde{u}_{j}(0)|=|Du_{j}(0)|\to \infty$, and
$\sup_{B_{h_{j}/(K\lambda_{j})}}\widetilde u_{j}= h_{j}/u_{j}(0)$ and $\widetilde u_{j}(0)=1$. Then if $h_{j}/u_{j}(0)$ is
bounded above by some constant $\beta<\infty$ we could deduce from Corollary~\ref{holdercor} that $|D
u_{j}(0)|(=\widetilde u_{j}(0)|D\widetilde u_{j}(0)|)$ is bounded above by a constant depending only on $m,n$, and
$\beta$, a contradiction since $|Du_{j}(0)|\to\infty$.  Thus~\ref{hoveru} is proved.

Now with $h_{j}$ as above, consider the new rescaling
\ab{3pt}{3pt}{%
$$%
\widetilde{G}_j=h_{j}^{-1} G_j.
$$}%
Observe that then $\widetilde{G}_{j}=\graph \widetilde{u}_{j}$, where $\widetilde u_{j}(x)=h_{j}^{-1}u(h_{j}x)$ on
$B_{1/K}(0)$, and
\ab{4pt}{-4pt}{%
$$%
\text{$\widetilde{u}_{j}(0)=u_{j}(0)/h_{j}\to 0$,\, 
$|D\widetilde{u}_{j}(0)|=|Du_{j}(0)|\to \infty$, \, $\sup_{B_{1/K}(0)}\widetilde u_{j}=1$}
$$}%
by~\ref{nonempty} and \ref{hoveru}.

\enlargethispage{0.3cm}

By Corollary~\ref{holdercor} the $\widetilde u_{j}$ are equicontinuous on $B_{\rho}(0)$ for each $\rho<1/K$ and hence
a subsequence converges locally uniformly in $B_{1/K}(0)$ to a singular solution $\widetilde u$ of \ref{sme} and the
graph $\widetilde G$ of $\widetilde u$ intersects the cone $W$ at each height $0< t< 1$, so every tangent cone of
$\widetilde G$ at $(0,0)$ has rays of slope $\ge K$, thus contradicting Theorem~\ref{tan-cones}, provided we choose
$K=\max\{M/(1-\theta),K_{0}+1\}$ with $K_{0}$ as in Lemma~\ref{bd-slope}.
\end{proof}


\section[dim\,sing\,$u\le n-2$] {dim\,sing\bm{$\,u\le n-2$}}
\label{dimsing} 

\begin{state}{\tl{sing0}\,Theorem.}%
  Suppose $u$ is a singular solution of the \emph{SME} (i.e.\ a singular solution of \emph{\ref{sme}}) in the domain
$\Omega\subset\R^{n}$. Then $\sing u(=\{x\in B_{1}:u(x)=0\})$ has Hausdorff dimension $\le n-2$; in fact, for each
closed ball $\overline{\!B}_{\rho}(y)\subset \Omega$, $\sing u\cap \overline{\!B}_{\rho}(y)$ can be written as a finite union
of locally compact (i.e.\ intersection of a compact with an open set) subsets, each of which has finite $(n-2)$-dimensional
Hausdorff measure in a neighborhood of each of its points.
\end{state} 

\tl{weak-check}\,{\bf{}Remark.} Using the bound on the singular set in the above theorem, together with the gradient
estimate of Theorem~\ref{grad0}, we can now check that singular solutions $u$ of~\ref{sme} on a domain $\Omega$
automatically have $1/u\in L^{1}_{\text{loc}}(\Omega)$ and are weak solutions of~\ref{sme}; i.e.~\ref{sme}$^{\prime}$
holds.  To check this we replace $\zeta$ in~\ref{sme}$^{\prime}$ by $\zeta\chi_{j}$ where $\chi_{j} \equiv 0 $ in the
$(1/j)$ neighborhood of support $\zeta\cap \{x\in \Omega:u(x)=0\}$, $\chi_{j} \equiv 1$ outside the $(2/j)$ neighborhood
of support $\zeta\cap \{x\in \Omega:u(x)=0\}$, and $\tint_{\R^{n}}|D\chi_{j}|\to 0$ as $j\to \infty$; such $\chi_{j}$ exist
because ${\cal{}H}^{n-1}(\sing u)=0$ by the above theorem.
 
\smallskip

\begin{proof}{\bf{}Proof of Theorem~\ref{sing0}:} For $K>0$ and $\Omega\subset\R^{n}$ open, let
  \ab{3pt}{3pt}{%
  $$%
  {\cal{}M}_{K,\Omega}= \{u: u\text{ is a singular solution of \ref{sme} on }\Omega \text{ with }\sup|Du|\le K\}
 $$}%
and
\ab{0pt}{3pt}{%
$$%
{\cal{}M}_{K}=\cup_{\text{open }\Omega\subset\R^{n}}{\cal{}M}_{K,\Omega},\quad S({\cal{}M}_{K})=\{S(u):u\in
{\cal{}M}_{K}\},
$$}%
where as usual $S(u)$ denotes the symmetric graph of $u$.

Then the gradient bound of Theorem~\ref{grad0} guarantees that each singular solution $u$ of the SME \ref{sme} on an
open $\Omega\subset\R^{n}$ must have $u|B_{\rho}(y)\in{\cal{}M}_{K,B_{\rho}(y)}$ for some $K$, provided $\overline
B_{\rho}(y)\subset\Omega$.  Also, with the aid of the Arzela-Ascoli lemma, one can readily check that, for each $K>0$,
$S({\cal{}M}_{K})$ is a ``multiplicity one class'' of stationary minimal hypersurfaces in $\R^{m+n}$ in the sense
of~\cite{Sim93}.  Then, as discussed in \cite{Sim93}, there is an integer $q\in \{0,\ldots,n+m-2\}$ such that $\dim\sing
M\le n+m-1-q$ for each $M\in S({\cal{}M}_{K})$, where $q$ is is the maximum integer such that there is
$\C=\C_{0}\times\R^{n+m-1-q}\in S({\cal{}M}_{K})$, where $\C_{0}$ is a minimal hypercone in $\R^{q+m}$ with
$\sing\C_{0}=\{0\}$ which is invariant under rotations of the last $m$ coordinates.  Then of course $q\neq 0$, and in
terms of the functions $u\in{\cal{}M}_{K}$ this says that $\dim\sing u\le n-q$ where $q$ is the minimum integer $\ge 1$
such that there is a homogeneous degree 1 singular Lipschitz solution $u$ of the SME \ref{sme} on $\R^{q}$.  $q$ is not
equal to $1$ because, by the discussion in \- \S\ref{intro}, there are no singular solutions of \ref{sme} in case $n=1$. On
the other hand as discussed in \S\ref{intro} there is the homogeneous degree 1 singular solution $\sqrt{m-1}|x|$ of
\ref{sme} in $\R^{2}$. So $q=2$ and hence each $u\in{\cal{}M}_{K}$ has $\sing u$ of Hausdorff dimension $\le n-2$ as
claimed.

The remaining rectifiability claims are true by~\cite{Sim95}.~\end{proof}

\providecommand{\bysame}{\leavevmode\hbox to3em{\hrulefill}\thinspace}


\begin{thebibliography}{BDG69}

\bibitem[BDG69]{BomDG69} E.~Bombieri, E.~De\hskip.05in{G}iorgi, and E.~Giusti, \emph{Minimal cones and the
{B}ernstein problem}, Invent. Math. \textbf{7} (1969), 243--268.

\bibitem[BG72]{BomG72} E.~Bombieri and E.~Giusti, \emph{Harnack's inequality for elliptic differential equations on
minimal surfaces}, Invent.\ Math. \textbf{15} (1972), 24--46.

\bibitem[CHS84]{CafHS84} L.~Caffarelli, R.~Hardt, and L.~Simon, \emph{Minimal surfaces with isolated singularities},
Manuscripta Math. \textbf{48} (1984), 1--18.

\bibitem[Dei85]{Dei85} K.~Deimling, \emph{Nonlinear functional analysis}, Springer-Verlag, 1985.

\bibitem[DH90]{DieH90} U.~Dierkes and G.~Huisken, \emph{The n-dimensional analogue of the catenary: existence and
nonexistence.}, Pacific J. Math. 141 (1990), no. 1, 47--54.

\bibitem[DH96]{DieH96} \bysame \emph{The $N$-dimensional analogue of the catenary}, Geometric analysis and the
calculus of variations, \textbf{12} (1996) Int. Press, Cambridge, MA.

\bibitem[EG]{EG} L.~Evans and R.~Gariepy, \emph{Measure Theory and Fine Properties of Functions}, (Revised Edition)
CRC Press, 2015.

\bibitem[GT]{GT} D.~Gilbarg and N.~Trudinger, \emph{Elliptic partial differential equations of second order}, 2nd ed.,
Springer, Berlin, 1983.

\bibitem[Il96]{Il96} T.\ Ilmanen, \emph{A strong maximum principle for singular minimal hypersurfaces}, Calc.  Var.  Partial
Differential Equations \textbf{4} (1996), 443--467.

\bibitem[KS89]{KorS89} N.\ Korevaar and L.~Simon, \emph{Continuity estimates for solutions to the prescribed curvature
{D}irichlet problem}, Math. Zeit. \textbf{197} (1989), 457--464.

\bibitem[Law72]{Law72} H.B. Lawson, \emph{The equivariant {P}lateau problem and interior regularity}, Trans.\ Amer.\
Math.\ Soc. \textbf{173} (1972), 231--249.

\bibitem[MicS73]{MicS73} J.H.~Michael and L.~Simon, \emph{Sobolev and Mean Value Inequalities on generalized
submanifolds of $\R^n$}, Comm. Pure Appl. Math \textbf{13} (1973)

\bibitem[Sis68]{Sis68} J.~Simons, \emph{Minimal varieties in {R}iemannian manifolds}, Ann.\ of Math.\ \textbf{88} (1968),
62--105.

\bibitem[Sim76]{Sim76} L.~Simon, \emph{Interior gradient bounds for non-uniformly elliptic equations}, Indiana Univ.\
Math.\ J. \textbf{25} (1976), 821--855.

\bibitem[Sim83]{Sim83} L.~Simon, \emph{Lectures on Geometric Measure Theory}, Proc. Centre Math. Anal.  Austral.  Nat.
 Univ. 3 (1983).

\bibitem[Sim87]{Sim87} L.~Simon, \emph{A Strict Maximum Principle For Area Minimizing Hypersurfaces}, J.~Differential
Geometry. \textbf{26} (1987), 327-335.

\bibitem[Sim93]{Sim93} L.~Simon, \emph{Cylindrical tangent cones and the singular set of minimal submanifolds},
J.~Differential Geometry. \textbf{38} (1993), 585-652.

\bibitem[Sim95]{Sim95} L.~Simon, \emph{Rectifiability of the singular sets of multiplicity 1 minimal surfaces and energy
minimizing maps}, Surveys Differential Geometry. \textbf{II} (1995), 246-305.

\bibitem[SS81]{SS81} R.~Schoen and L.~Simon, \emph{Regularity of stable minimal hypersurfaces}, Comm. Pure and Appl. 
Math.\ \textbf{34}, (1981), 742--797.

\bibitem[Wic14]{Wic14} N.~Wickramasekera, \emph{A General Regularity Theory for Stable Codimension~\emph{1} Integral
Varifolds}, Annals of Mathematics \textbf{179} (2014), 843--1007.
 
\end{thebibliography}
\end{document}